\documentclass[10pt]{article}
\usepackage{amsfonts}
\usepackage{latexsym,amssymb}
\usepackage{amsmath, amsbsy}
\usepackage{amsopn, amstext}
\usepackage{graphicx, color, epstopdf}
\usepackage{threeparttable}
\usepackage{hyperref}

\usepackage{array}
\usepackage{bm}
\usepackage{multirow}
\usepackage{booktabs}

\usepackage{subfigure}
\usepackage{tabulary}

\def\R{\mathbb{R}}
\newtheorem{corollary}{Corollary}

\newtheorem{remark}{Remark}

\newtheorem{proposition}{Proposition}

\numberwithin{equation}{section}
 \numberwithin{Lem}{section}
 \numberwithin{Defi}{section}
 \numberwithin{Theo}{section}
 \numberwithin{Rem}{section}
 \numberwithin{Coro}{section}
 \numberwithin{Fig}{section}

\newcommand{\norm}[1]{\left\Vert#1\right\Vert}

\newcommand{\abs}[1]{\left\vert#1\right\vert}

\newcommand{\ds}{\,\mathrm{d}s}
\newcommand{\dx}{\,\mathrm{d}x}
\renewcommand{\i}{\mathrm{i}}

\title{Numerical solution of a one-dimensional nonlocal Helmholtz equation by Perfectly Matched Layers\thanks{This work is supported in NSFC under grants No. 11771035, 12071401 and NSAF U1930402, Natural Science Foundation of Hunan Province No. 2019JJ50572, Natural Science Foundation of Hubei Province No. 2019CFA007 and Xiangtan University 2018ICIP01. }}

\author{
 Yu Du\thanks{Department of Mathematics, Xiangtan University, Hunan, 411105, China({\tt duyu@xtu.edu.cn})}
\and Jiwei Zhang\thanks{School of Mathematics and Statistics, and Hubei Key Laboratory of Computational Science, Wuhan University, Wuhan 430072, China.
({\tt jiweizhang@whu.edu.cn})}
}
\date{}

\begin{document}

\maketitle

\begin{abstract}


We consider the computation of a nonlocal Helmholtz equation by using Perfectly Matched Layer (PML). We first derive the nonlocal PML equation by extending PML modifications from the local operator to the nonlocal operator of integral form. We then give stability estimates of some weighted average value of the nonlocal Helmholtz solution and prove that (i) the weighted average value of the nonlocal PML solution decays exponentially in PML layers in one case; (ii) in the other case, the weighted average value of the nonlocal Helmholtz solution itself decays exponentially outside some domain. Particularly for a typical kernel function $\gamma_1(s)=\frac12 e^{-| s|}$, we obtain the Green's function of the nonlocal Helmholtz equation, and use the Green's function to further prove that (i) the nonlocal PML solution decays exponentially in PML layers in one case; (ii) in the other case, the nonlocal Helmholtz solution itself decays exponentially outside some domain. Based on our theoretical analysis, the truncated nonlocal problems are discussed and an asymptotic compatibility scheme is also introduced to solve the resulting truncated problems. Finally, numerical examples are provided to verify the effectiveness and validation of our nonlocal PML strategy and theoretical findings.



\end{abstract}

\section{Introduction}\label{sec_intro}
The development of nonlocal models has grown impressively over the last decade because of its huge potential of emerging applications in various research areas, such as the peridynamical theory of continuum mechanics, the nonlocal wave propagation, and the modeling of nonlocal diffusion process \cite{bobaru2010the,du2012analysis,silling2000reformulation,weckner2005the,zhou2010mathematical}. In this paper, we consider the computation of a nonlocal Helmholtz equation on the whole real axis 
\begin{align}
\mathcal{L}_\delta u(x) - k^2u(x) = f(x), \quad x\in\R, \label{eq:nonlocalHelmholtz}
\end{align}
where $k$ is a constant related to the traditional \emph{wavenumber} for the local Helmholtz equation, the source $f\in L^2(\R)$ is supported on $\Omega:=(-l,l)$, and the nonlocal operator $\mathcal L_\delta$ is defined as
\begin{align}
	\mathcal{L}_\delta u(x) = \int_\R \big( u(x)-u(y) \big) \gamma_\delta(y-x) \mathrm{d}y.
	\label{eq:nonlocalOperator}
\end{align}
The kernel function $\gamma_\delta$ in \eqref{eq:nonlocalOperator} is determined by a rescaling of a parent kernel $\gamma_1$ through
\begin{equation}
\gamma_\delta( s)=\frac{1}{\delta^3}\gamma_1\big(\frac{ s}{\delta}\big),
\end{equation} 
where the parameter $\delta$ represents the range/radius of nonlocal interaction, and $\gamma_1( s)\in L^1{(\mathbb{R})}$ is piecewisely smooth, and satisfies:
\begin{itemize}
\item nonnegativeness: $\gamma_1( s)\geq0$;
\item symmetry in $ s$: $\gamma_1( s)=\gamma_1(- s)$;
\item finite horizon: $\exists\; l_\gamma>0$, such that $\gamma_1( s)=0$ if $| s|>l_\gamma$;
\item the second moment condition $\frac12 \int_\R s^2\gamma_1( s)d s=1.$
\end{itemize}


Recently, much works are carried out for the simulation of nonlocal problems with free or fixed boundary conditions. There are applications in which the simulation of an infinite medium may be useful, such as wave or crack propagation in whole space. 
The nonlocal Helmholtz equation can be used to describe the nonlocal wave propagation. In fact, it can be derived from the nonlocal wave equation
\begin{align}
	\left( \partial_t^2 + \mathcal{L}_\delta \right) u(x,t) = f(x,t),\quad x\in\R, \label{eq:nonlolocalWaveEquation}
\end{align}
where $u(x,t)$ represents the displacement field, $f(x,t)$ is the space-time source term with compact support at all time. If we make the ansatz that $f(x,t)$ is a superposition of the time-harmonic sources $f(x)e^{-\i k t}$. Then, for each $k$, the corresponding mode $u(x)e^{-\i kt}$ satisfies
\begin{align*}
	\mathcal{L}_\delta \big(u(x)e^{-\i kt}\big) - k^2 u(x)e^{-\i kt} = f(x)e^{-\i k t}.
\end{align*}
Thus the time domains solution $u(x,t)$ is the sum of the time-harmonic modes $u(x)e^{-\i kt}$ over all possible values of $k$.

The aim of the paper is to develop an efficient computation of a nonlocal Helmholtz equation on the whole real axis. Absorbing boundary conditions (ABCs) are a successful approach to simulate the wave behaviors of a physical domain of interest by imposing a suitable boundary condition, to absorb the impinging wave at artificial boundaries. For the construction of tractional ABCs, it is well studied for local problems \cite{AlpertGreengardHagstrom, HanHuang, hagstrom2008high, grote1995exact, teng2003exact, houdehan2003exact, arnold2005approximation}, and there are also much progress for nonlocal problems, see \cite{DuHanZhangZheng, ZYZD16, ZhengHuDuZhang, wildman2012a}. In this paper, we consider to construct the perfectly matched layer (PML) as ABCs for a 1D nonlocal Helmholtz equation. The PML, originally proposed by Berenger \cite{Berenger1994}, has the two important features: (i) the wave in a special designed layer decays exponentially, and (ii), if the wave reflects off the truncated boundary, the returning waves after one round trip through the absorbing layer are very tiny \cite{AL,cw,berenger1996three,collino1998the,turkel1998absorbing,ls98,cw03,Chen2005,bw05,li2018fem,bramble2006analysis,hohage2003solving}.
Specifically for peridynamics, Wildman and Gazonas \cite{wildman2012a} present a PML by treating the nonlocal kernel as the convolution of the displacement with the second derivative of a nascent Dirac delta distribution. In this paper, we derive a PML different from \cite{wildman2012a} and has a simpler structure, and more importantly, give the stability estimate theoretically.

 The contributions of this paper are given as follows:
\begin{itemize}
\item We extend the strategy of the PML method from the local operator to nonlocal operator, and obtain a nonlocal PML equation directly from the weak form of nonlocal Helmholtz equation. The resulting nonlocal PML equation converges to the corresponding local PML equation while the nonlocal interaction horizon vanishes. Such consistency is useful to demonstrate the validation/verification of our analytical continuation for the nonlocal operator.
\item The properties of the nonlocal Helmholtz solution and nonlocal PML solution are analyzed. To do so, we introduce the weighted average value $u^w$ with respect to the nonlocal Helmholtz solution $u$, defined by
\begin{align*}
	u^w(x) := \int_{\R} u(t+x)\kappa(t)\mathrm dt, 
\end{align*}
where the weight $\kappa(t)$ is given in \eqref{w1}-\eqref{w2}. For general kernel function $
\gamma_\delta$, the stability estimates of $u^w$ and its analytical continuation $\tilde{u}^w$ are established. Specifically, we prove in Corollary~\ref{th:barq}: (1) the weighted average value $\tilde{u}^w$ of the nonlocal PML solution exponentially decays in PML layers under a quantitative condition depending on the \emph{wavenumber} $k$ (taken smaller value) and the kernel $\gamma_\delta$; (2) under another quantitative condition, the weighted average value $u^w$ itself of the nonlocal Helmholtz solution exponentially decays outside of a domain large enough. The  result (1) shows the PML is efficient to absorb the waves when they impinge the PML layers. The result (2) suggests that one can directly truncate the nonlocal Helmholtz equation by putting forced boundary layers (i.e. homogeneous Dirichlet boundary constrains) at sufficiently large $x$ instead of using PML. 

\item In particular, the refined estimates of the nonlocal Helmholtz solution $u$ and its PML solution $\tilde{u}$ are further established for a typical kernel $\gamma_1=\frac12 e^{-|s|}$. To this end, we first derive the exact formula of the Green's function for nonlocal Helmholtz equation~\eqref{eq:nonlocalHelmholtz}, and then present the stability estimates of $u$ and $\tilde{u}$ instead of their weighted average values. Specifically, we prove in Corollary~\ref{th:expnonlocalPML} that: (1) if $\delta k<1$, the nonlocal PML solution exponentially decays in PML layers; (2) if $\delta k>1$, the nonlocal Helmholtz solution $u$ itself exponentially decays outside $\Omega$. The resulting Green's function for this typical kernel shows an attractive property that the local and nonlocal Helmholtz solutions have an intimate connection, i.e., they can be expressed by each other. 
\end{itemize}

The behaviors of nonlocal Helmholtz (PML) solutions and Green's function for the Helmholtz equation with $\gamma_1=\frac12 e^{-|s|}$ are newly obtained, although the Green's functions of some nonlocal problems have been previously studied \cite{wang2017static, hanson2008dyadic, weckner2005the, mikata2012analytical, weckner2009green} which are mainly based on Fourier transforms and Laplace transforms, but can not be extended to the nonlocal Helmholtz equation.

On the other hand, the nonlocal operator $\mathcal{L}_\delta$ has an intimate connection with some local differential operator. More specifically, as the nonlocal horizon vanishes, see~\cite{du2012analysis, zhou2010mathematical, du2019nonlocal, du2018nonlocal}, the nonlocal operator $\mathcal{L}_\delta$ converges to the second-order differential operator 
\begin{equation} \label{eq:NTL}
\lim_{\delta\rightarrow 0}\mathcal{L_\delta} u(x) = -\partial_x^2u(x) :=\mathcal{L}_0 u(x)
\end{equation}
under the above second moment condition.  Such consistency is quite useful not only for the modeling, but also for the validation/verification of numerical simulations. For discrete schemes, it is useful to make use of the asymptotic compatibility (AC) scheme, a concept developed in \cite{TianDu, TianDu2} and further extended in \cite{du2019a,TianDu3}, to discretize the nonlocal operator to preserve such (analogous) limit \eqref{eq:NTL} in a discrete level. Recently, Du et al. \cite{du2018nonlocal, DuHanZhangZheng} proposed an AC scheme to discretize the nonlocal operator in unbounded multi-scale mediums. In this paper we extend the method \cite{du2018nonlocal} to discretize the nonlocal operator with PML modifications.

The outline of this paper is as follows. In Section~\ref{sec:PMLforHelmholtz} we propose our PML technique for the nonlocal Helmholtz equation by extending complex coordinate transforms for the local Helmholtz equation. Section~\ref{sec:exponentiallyDecayingWaves} is devoted to analyzing the nonlocal Helmholtz solution and give the truncated nonlocal problems based on our theoretical analysis. In particular, for a special kernel $\gamma_1( s)=\frac12 e^{-| s|}$ we give the Green's function of the nonlocal problem. In Section~\ref{sec:Discretization}, we introduce the AC scheme for numerically solving the truncated problems. In Section~\ref{sec:numericalTests}, some numerical tests are provided to verify the theoretical results and the effectiveness of our PML.

\section{PML for the local and nonlocal Helmholtz equations}\label{sec:PMLforHelmholtz}
In this section we first review the classical PML method for local Helmholtz equations, and recall some properties of local PML solutions. We then employ the weak form of nonlocal equation~\eqref{eq:nonlocalHelmholtz} to derive the corresponding nonlocal PML equation, and finally consider the local limits between the local and nonlocal PML equations as the nonlocal horizon $\delta$ vanishes.

\subsection{PML for the local Helmholtz equation}
We now recall the PML method for a local Helmholtz equation
\begin{align} 
	\mathcal{L}_0 u_{loc}(x) - k^2 u_{loc}(x) = f(x) \label{eq:localHelmholtz}
\end{align}
with Sommerfeld radiation condition. As discussed in \cite{cw, cjm97, bramble2006analysis, collino1998the, ls98}, the PML modifications can be viewed as a complex coordinate stretching of the original scattering problem by constructing an analytic continuation to the complex plane, namely, 
\begin{align}
	\tilde x(x) = \int_{0}^x \omega(t)\mathrm{d}t = \int_{0}^x 1+\i \sigma(t)\mathrm{d}t, \label{eq:omegax}
\end{align}
where $\sigma$, called the absorption function, is equal to $0$ in $\Omega$ and positive outside $\Omega$. $\i=\sqrt{-1}$ is the imaginary unit.

Setting $\tilde u_{loc}(x)=u_{loc}(\tilde x)$, and using the relation 
\begin{align*}
\partial_x\rightarrow \frac{1}{\omega(x)} \partial_x = \frac{1}{1+\i\sigma(x)} \partial_x,
\end{align*}
the modified modal solution $\tilde u_{loc}(x)$ satisfies the governing equation 
\begin{align}\label{eq:localPML} 
	\tilde{\mathcal L}_0 u_{loc}(x)-k^2\omega(x) u_{loc}(x)= f(x)
\end{align}
with the local PML operator given by
\[
	\tilde{\mathcal L}_0 u_{loc}(x)=-\partial_x\left( \frac{1}{\omega(x)} \partial_x\tilde u_{loc}(x)\right).
\]


The PML method for local Helmholtz equations has been well studied, see \cite{li2018fem,ms10}. 
Here we recall some properties of solutions as follows:
\begin{itemize}
\item[(i).] The local Helmholtz solution $u_{loc}$ is oscillating and possibly doesn't vanish as $|x|\to\infty$.
\item[(ii).] The local PML solution $\tilde u_{loc}(x)$ is an analytic continuation of $u_{loc}(x)$ in the complex coordinate. The solution $\tilde u_{loc}(x)$ is unique in the sense that $\tilde u_{loc}(x)=u_{loc}(x)$ over the domain $\Omega$.
\item[(iii).] The stability estimate of the solution to the local Helmholtz equation is given as 
\begin{align*}
\|u_{loc}\|_{H^1(\Omega)} + k\|u_{loc}\|_{L^2(\Omega)} \leq C \|f\|_{L^2(\Omega)}.
\end{align*}
\item[(iv).] The analytical continuation changes oscillating waves into exponentially decaying waves outside the region of interest,
\begin{align}
|\tilde u(x)| \leq C e^{-k|\int_0^x\sigma(t)\mathrm dt|}\|f\|_{L^2(\Omega)}\quad \mathrm{for}\ |x|>l. \label{eq:localPMLexpDecay}
\end{align}
\end{itemize}

\subsection{PML for the nonlocal Helmholtz equation}
The original PML technique is proposed for the standard PDEs, such as the local Helmholtz equation shown above and the Maxwell's equation. We point out that PML is for the problems with some far field boundary conditions. An typical example is the Sommerfeld radiation condition for the ``local'' Helmholtz equation in a homogeneous medium. However, it is unknown for the nonlocal Helmholtz equation~\eqref{eq:nonlocalHelmholtz}. Therefore, in this paper we simply assume that a suitable boundary condition at infinity is imposed to exclude energy incoming from infinity and only to allow energy outgoing to infinity. Based on this assumption, we here manage to extend the PML to solve the nonlocal Helmholtz equation.

To do so, we introduce the PML by considering the weak form of Eq.~\eqref{eq:nonlocalHelmholtz}
\begin{align}
	(\mathcal{L}_\delta u,v) - k^2(u,v) = (f,v),\quad \forall v\in C_0^\infty(\R), \label{eq:weakform}
\end{align}
where $(\cdot,\cdot)$ denotes the inner product in the complex valued $L^2$-space, and 
\begin{align}
	(\mathcal{L}_\delta u,v) = \frac12 \int_\R \int_\R \big( u(x)-u(y) \big) \big( \bar v(x)-\bar v(y) \big) \gamma_\delta(y-x) \mathrm{d}y\mathrm{d}x. \label{eq:weakcall}
\end{align}
Here $\bar{v}(x)$ represents the complex conjugate of $v(x)$. We now apply the same transform as \eqref{eq:omegax} and immediately produce the corresponding differential forms by 
\begin{align}
	&x\rightarrow \tilde x=\int_0^x \omega(t)\mathrm{d}t, &&y\rightarrow \tilde y=\int_0^y \omega(t)\mathrm{d}t, \label{eq:complexStreching} \\
	&\mathrm{d}x \rightarrow \omega(x)\mathrm{d}x, &&\mathrm{d}y\rightarrow \omega(y)\mathrm{d}y.
\end{align}
Setting $\tilde u(x)=u(\tilde x)$, $\tilde u(y)=u(\tilde y)$, $\tilde{\bar v}(x)=\bar v(\tilde x)$ and $\tilde{\bar v}(y)=\bar v(\tilde y)$, we can transform~\eqref{eq:weakform} into the following nonlocal equation with PML modifications
\begin{align}
	 & \frac12 \int_\R \int_\R \big( \tilde u(x)-\tilde u(y) \big) \big( \tilde{\bar v}(x)-\tilde{\bar v}(y) \big) \gamma_\delta(\tilde y-\tilde x) \omega(x)\omega(y)\mathrm{d}y\mathrm{d}x \label{eq:PML} \\
	 &\quad\quad\quad\quad\quad\quad\quad\quad\quad\quad\quad\quad- k^2\int_\R \tilde u(x)\tilde{\bar v}(x)\omega(x)\mathrm{d}x = \int_\R f(x)\tilde{\bar v}(x)\omega(x)\mathrm{d}x. \notag
\end{align}
Using the facts $\omega(x)=1 \; \forall x\in \Omega$ and $\mathrm{supp} f(x)\subset \Omega$, we have the strong form of \eqref{eq:PML} as 
\begin{align} 
	\tilde{\mathcal{L}}_\delta\tilde u(x) - k^2\omega(x)\tilde u(x) = f(x), \label{eq:PMLeq}
\end{align}
where the nonlocal operator with PML modifications is given by
\begin{align}
	\tilde{\mathcal{L}}_\delta \tilde u(x) = \int_\R \big( \tilde u(x)-\tilde u(y) \big) \gamma_\delta(\tilde y-\tilde x) \omega(x)\omega(y) \mathrm{d}y. \label{eq:nonlocalPMLoperator}
\end{align}
Thus, we have that the solution $\tilde u(x)$ of Eq.~\eqref{eq:PMLeq} is an analytic continuation of the solution $u(x)$ of Eq.~\eqref{eq:nonlocalHelmholtz} in the complex coordinate, and it holds that $\tilde u(x)=u(x)$ for $x\in\Omega$. We emphasize that the kernel function $\gamma_\delta$ in the nonlocal PML operator~\eqref{eq:nonlocalPMLoperator} must be the analytic continuation of the original kernel in complex coordinates. From Figure~\ref{fig:ex1_solution_diffkernel} in the section of numerical examples, one can see that such analytic continuation of the kernel plays an important role to get the correct solution in simulations. We give the analytic continuation of two typical kernels in Section~\ref{sec:numericalTests}, and the analytic continuation of general kernels still needs to be studied.

On the other hand, it is well known that the nonlocal operator $\mathcal{L}_\delta$ has an intimate connection with the local differential operator shown in Eq.~\eqref{eq:NTL}. It is interesting to ask if the nonlocal PML equation~\eqref{eq:PMLeq} converges to its corresponding local PML equation \eqref{eq:localPML} as the nonlocal interaction horizon $\delta$ vanishes. Such consistency is useful to demonstrate the validation/verification of our analytical continuation for the nonlocal operator. In fact, by using the second moment condition of the kernel function and Taylor expansion, it is straightforward to verify that 
\begin{align}\label{eq:nPlimit}
	\lim_{\delta\rightarrow0^+} (\tilde{\mathcal{L}}_\delta\tilde u(x),v(x)) = (\tilde{\mathcal L}_0 \tilde u(x),v(x)).
\end{align}
Thus, the local limit \eqref{eq:nPlimit} implies that the nonlocal PML equation \eqref{eq:PMLeq} will converge to the corresponding local PML equation \eqref{eq:localPML} as the nonlocal interaction horizon tends to zero. 

\section{The exponentially decaying waves} \label{sec:exponentiallyDecayingWaves}
In the previous section, the ``a priori" estimates of solutions are introduced for the local Helmholtz equation and its PML equation, see \cite{li2018fem,ms10}. These properties are useful for the further theoretical analysis of the PML method. Here we also want to study the solution properties of the nonlocal Helmholtz equation \eqref{eq:nonlocalHelmholtz} and the nonlocal PML equation \eqref{eq:PML}. Generally, the properties of nonlocal solutions are not trivial to be explored since it is hard to exactly express the Green's function for the general kernel.
%
To do so, we first introduce the following function
\begin{align}
	G_{x_0}(x) = \begin{cases}
		C_1(x_0) e^{-\i\tilde kx}, & \quad x\leq x_0 \quad \text{with} \quad C_1(x_0)=- \frac{e^{\i \tilde k x_0}}{2\i \tilde k}, \\
		C_2(x_0) e^{\i\tilde kx}, & \quad x>x_0\quad \text{with} \quad C_2(x_0)=- \frac{e^{-\i \tilde k x_0}}{2\i \tilde k}.
	\end{cases}
	\label{eq:green}
\end{align}
Here $\tilde k$ is the solution of the equation 
\begin{align} \label{eq:relation}
	\mu(\tilde k)=k^2,\qquad \mathrm{where}\ \mu(\tilde k)&:= \frac{1}{\delta^2} \int_\R \left(1-e^{\i\tilde k\delta s}\right)\gamma_1(s)\ds.
\end{align}

The identity \eqref{eq:relation} implies that $\tilde k$ is only dependent of the kernel $\gamma_1$, $\delta$ and the \emph{wavenumber} $k$, and $\tilde k$ converges to $k$ as $\delta \to0$ under the assumption of second moment condition of the kernel $\gamma_1$.

In the following sections, we assume that there exists the solution $\tilde k$ for the equation $\mu(\tilde k)=k^2$, which still need to be studied for general kernels and \emph{wavenumber} $k$. Under this assumption, we know that it also holds $\mu(-\tilde k)=k^2$ since the kernel $\gamma_1$ is symmetric, which implies that there exists $\tilde k$ in the set of solutions to Eq.~\eqref{eq:relation} such that $\tilde k\in\R^+$ or $\Im(\tilde k)>0$ if $k^2\neq \mu(0)$. Specifically, there exists a positive value $k_0$ which is
$$k_0= \sup_{\tilde k\in \R^+}\sqrt{\mu(\tilde k)},$$
such that $\tilde k\in\R^+$ if $k<k_0$, and $\Im(\tilde k)>0$ if $k>k_0$. In this paper, we only need to take the root $\tilde k$ of identity \eqref{eq:relation} as a positive number or a complex number with positive imaginary part.

To make clear the relation between $k$ and $\tilde k$, we here show two examples.
\begin{itemize}
 \item[] {\bf Example 1.} Taking the kernel $\gamma_1( s)=\frac{1}{2} e^{-| s|}$, we have 
	  \begin{align*}
		  \mu(\tilde k)= \tilde k^2/(\delta^2\tilde{k}^2+1). 
	  \end{align*}
	  The direct calculation shows that one typical solution of Eq.~\eqref{eq:relation} is
	  \begin{align}
	  	\tilde{k}=k\sqrt{\frac{1}{1-(\delta k)^2}}. \label{eq:ktildekdelta}
	  \end{align}
	  It's clear that $\tilde k\in\R^+$ if $k<k_0=1/\delta$ and $\tilde k=|\tilde k|\i$ if $k>k_0= 1/\delta$.
	  \item[] {\bf Example 2.} Taking the kernel $\gamma_1( s) = \frac{4}{\sqrt{\pi}} e^{- s^2}$, we have 
 \begin{align*}
\mu(\tilde k) = \frac{4}{\delta^2}\left(1-e^{-(\delta\tilde k)^2/4}\right).
\end{align*}
One typical solution of Eq.~\eqref{eq:relation} is
\begin{align}
\tilde k=\frac{2}{\delta} \sqrt{-\log\left(1-(k\delta)^2/4\right)}.
\label{eq:gausskerneltildek}
\end{align}
It is clear that $\tilde k\in\R^+$ if $k<k_0=2/\delta$, and $\Im(\tilde k)>0$ if $k>k_0=2/\delta$.
\end{itemize}

To explore the properties of the solution $u(x)$ to the nonlocal Helmholtz equation~\eqref{eq:nonlocalHelmholtz}, we now introduce the definition of a weighted average value with respect to $u(x)$ by
\begin{align}
	u^w(x) := \int_{\R} u(t+x)\kappa(t)\mathrm dt, \label{eq:barqdef}
\end{align}
where the weight is given by $\kappa(t) = w_1(t)$ for $t<0$ and $\kappa(t) = w_2(t)$ for $t>0$ with 
\begin{align}\label{w1}
&w_1(t) = \frac{1}{2\delta^2\i\tilde k}\int_{-\infty}^{\frac{t}{\delta}} \left(e^{\i\tilde k(t - \delta s)} - e^{-\i\tilde k(t-\delta s)} \right) \gamma_1(s) \ds = \frac{1}{\delta^2\tilde k}\int_{-\infty}^{\frac{t}{\delta}} \sin \tilde k(t - \delta s) \gamma_1(s) \ds, \\
&w_2(t) = \frac{1}{2\delta^2\i\tilde k} \int_{\frac{t}{\delta}}^{+\infty} \left( e^{-\i\tilde k (t-\delta s)} - e^{\i\tilde k(t-\delta s)} \right) \gamma_1(s) \ds = -\frac{1}{\delta^2\tilde k} \int_{\frac{t}{\delta}}^{+\infty} \sin \tilde k (t-\delta s) \gamma_1(s) \ds. \label{w2}
\end{align}

The weighted average value $u^w(x)$ satisfies the following identity. 
\begin{proposition}\label{Le1}
The weighted average value $u^w(x)$ defined by \eqref{eq:barqdef} satisfies
\begin{align}
u^w(x) = \int_{\R} G_{x}(y)f(y)\mathrm dy. \label{eq:baruExp}
\end{align}

\end{proposition}

For brevity, we leave the proof of Proposition \ref{Le1} in the appendix~\ref{Le1}.

\begin{corollary} \label{th:barq} 
Assume that $C$ is a constant only depending on $\Omega= (-l,l)$, and $\tilde{k}$ is the root of relation \eqref{eq:relation}. Then the weighted average value $u^w(x)$ of the solution $u(x)$ to the nonlocal Helmholtz equation \eqref{eq:nonlocalHelmholtz}, and its analytic continuation $\tilde{u}^w(x):=u^w(\tilde x)$ by complex coordinate transform \eqref{eq:complexStreching} have the following estimates:
	\begin{itemize}
		\item[(1)] If $\tilde k$ is positive real and $f\in L^2(\R)$, it holds that 
		\begin{align}
			&|u^w(x)|_{H^1(\Omega)} + \tilde k\|u^w(x)\|_{L^2(\Omega)} \leq C \|f\|_{L^2(\Omega)}, \label{eq:nonlocalAvgStaCase1}\\
			&|\tilde u^w(x)|\leq C e^{-\tilde k|\int_0^x \sigma(t)\mathrm dt|} \|f\|_{L^2(\Omega)},\quad \mathrm{for}\ |x|>l. \label{eq:nonlocalAvgExpDecayCase1}
		\end{align}
		\item[(2)] If $\tilde k$ is a complex number with $\Im(\tilde k)>0$ and $f\in L^2(\R)$, it holds that 
		\begin{align}
			&|u^w(x)|_{H^1(\Omega)} + |\tilde k|\|u^w(x)\|_{L^2(\Omega)} \leq C\|f\|_{L^2(\Omega)}, \label{eq:nonlocalAvgExpDecayCase2} \\
		&|u^w(x)| \leq C \frac{1}{|\tilde k|} e^{-\Im(\tilde k)(|x|-l)} \|f\|_{L^2(\Omega)},\quad \mathrm{for}\ |x|>l.
	\end{align}
	\end{itemize}
\end{corollary}

\begin{proof}
\emph{(1)} 
If $\tilde k$ is positive, \eqref{eq:green} is actually the Green's function of the following local Helmholtz equation \cite{ms10}
\begin{align*}
	-\partial_x^2\bar u(x) - \tilde k^2 \bar u(x) =      & f(x)\quad & & x\in\R,        \\
	\Big|\partial_r u^w -\i\tilde k u^w\Big|\to & 0  & & \mathrm{as}\ r=|x|\rightarrow+\infty.
\end{align*}
From Eq.~\eqref{eq:baruExp}, $u^w(x)$ is the solution of the Helmholtz problem above. By using the results for local problems in \cite{li2018fem,ms10} (see Eq.~\eqref{eq:localPMLexpDecay}), we can directly prove the stability estimate~\eqref{eq:nonlocalAvgStaCase1} for $u^w$ and the exponentially decaying estimate \eqref{eq:nonlocalAvgExpDecayCase1} for $\tilde u^w$.

\emph{(2)} From Eq.~\eqref{eq:green} and \eqref{eq:baruExp}, we get
	\begin{align*}
		\|u^w\|_{L^2(\Omega)}^2 =& \int_\Omega \Big| \int_\Omega G_x(y)f(y)\mathrm{d}y \Big|^2 \mathrm{d}x\leq \|f\|_{L^2(\Omega)}^2 \int_\Omega \int_\Omega \big|G_x(y)\big|^2\mathrm{d}y \mathrm{d}x \\
		\leq& \frac{1}{4|\tilde k|^2}\|f\|_{L^2(\Omega)}^2 \Big( \int_\Omega e^{-2\Im(\tilde k)x} \mathrm{d}x\int_{-l}^x e^{2\Im(\tilde k)y}\mathrm{d}y + \int_\Omega e^{2\Im(\tilde k)x} \mathrm{d}x\int_{x}^l e^{-2\Im(\tilde k)y}\mathrm{d}y \Big)\\
		=& \frac{1}{4|\tilde k|^2}\|f\|_{L^2(\Omega)}^2 \cdot \Big( \frac{2l}{\Im(\tilde k)}+\frac{1}{2\Im(\tilde k)^2}\big( e^{-4\Im(\tilde k)l}-1\big)\Big) \\
		=& \frac{2l^2}{|\tilde k|^2}\cdot \frac{ 4\Im(\tilde k)l+ e^{-4\Im(\tilde k)l}-1 }{(4\Im(\tilde k)l)^2} \|f\|_{L^2(\Omega)}^2\\
		\leq & \frac{l^2}{|\tilde k|^2} \|f\|_{L^2(\Omega)}^2,
	\end{align*}
where we have used the fact that $\sup_{t>0}\frac{t+e^{-t}-1}{t^2}=\frac12$.
	
By similar arguments, we have
	\begin{align*}
		|u^w|_{H^1(\Omega)}^2 \leq l^2 \|f\|_{L^2(\Omega)}^2.
	\end{align*}
We now prove that $u^w$ is exponentially decaying as $|x|\to\infty$. For $x>l$, it holds	
	\begin{align*}
		|u^w(x)| = & \left|\int_\Omega G_x(y)f(y) \mathrm{d}y\right| = \left|\int_{\Omega} C_0 e^{\i\tilde k(x-y)} f(y) \mathrm{d}y\right| \\
		\leq     & \frac{1}{2|\tilde k|} \left(\int_\Omega e^{-2\Im(\tilde k)(x-y)} \mathrm{d}y \right)^\frac12 \|f\|_{L^2(\Omega)}\\
		\leq& \frac{\sqrt{2l}}{2|\tilde k|} e^{-\Im(\tilde k)(x-l)} \|f\|_{L^2(\Omega)}. 
	\end{align*}
	Similarly, we can get $|u^w(x)| \leq \frac{\sqrt{2l}}{2|\tilde k|} e^{-\Im(\tilde k)(-x-l)} \|f\|_{L^2(\Omega)}$ for $x<-l$. This completes the proof.
\end{proof}

\begin{remark} For the main results in Corollary \ref{th:barq}, it is interesting to point out that 
	\begin{itemize}
	\item For the case of $\tilde k=|\tilde k|$, the estimate \eqref{eq:nonlocalAvgExpDecayCase1} shows that the analytic continuation of the weighted average value $\tilde{u}^w(x)$ decays exponentially as $|x|\rightarrow+\infty$. It maybe implies that the analytic continuation of the nonlocal Helmholtz solution $u(x)$ decays exponentially in some sense. 

	\item If $\Im(\tilde k)>0$, the weighted average value of $u(x)$ without any modification decays exponentially as $|x|\rightarrow+\infty$. It shows that the nonlocal Helmholtz solution $u(x)$ itself decays exponentially.

	\item These results are still valid for the case $l_\gamma=+\infty$, i.e., $\mathrm{supp} \;\gamma_1(s)\subset (-l_\gamma,l_\gamma)$. This means that $\gamma( s)$ decays as $| s|\rightarrow\infty$. In this situation, we can truncate the kernel when $\gamma( s)<\mathrm{\epsilon}$ for $| s|>\hat l_\gamma$ (a finite constant) for a given tolerance $\epsilon \; (e.g., \mathrm{\epsilon}=10^{-16})$. Taking the weighted average value with a truncated kernel by 
	\begin{align*}
		\hat{u}^w_\delta(x) := \int_{-\delta \hat l_\gamma}^{\delta \hat l_\gamma} u(t+x) \kappa(t) \mathrm{d}t,
	\end{align*}
we then have 
$$|\tilde{u}^w_\delta(x)-\hat{u}^w_{\delta}(\tilde x)|\leq 2\mathrm{\epsilon}\|\tilde u(x)\|_{L^1(\R)}/(\tilde k^2\delta^3).$$
This implies we can use the truncated kernel to replace the original one with a tolerance error.
\end{itemize}
\end{remark}


\subsection{A refined estimate for a typical exponential kernel $\gamma_1( s)=\frac12 e^{-| s|}$}
In the previous subsection, we considered the general compactly supported kernel function and have proved that the weighted average value of the nonlocal PML solution $\tilde u(x)$ or the nonlocal Helmholtz solution $u(x)$ decays exponentially as $x\to\infty$ (see Corollary~\ref{th:barq}). Here we consider a special kernel $\gamma_1( s)=\frac{1}{2} e^{-| s|}$, and directly analyze the behavior of $\tilde u(x)$ or $u(x)$ as $x\to\infty$ instead of their weighted average values. To do so, we first consider the Green's function $G_{x_0}^{\mathcal{L}_\delta}(x)$ for the following equation 
\begin{align}
	\mathcal{L}_\delta G_{x_0}^{\mathcal{L}_\delta}(x) -k^2 G_{x_0}^{\mathcal{L}_\delta}(x) = \mathcal{D}(x-x_0), \label{eq:nonlocalgreen}
\end{align}
where we use $\mathcal{D}(x-x_0)$ to represent the Dirac delta function differing from the nonlocal horizon $\delta$. The following corollary will present the explicit formula of $G_{x_0}^{\mathcal{L}_\delta}(x)$. 
\begin{corollary} \label{th:greenexpgamma} 
	For the kernel function $\gamma_1( s)=\frac{1}{2} e^{-| s|}$, the Green's function $G_{x_0}^{\mathcal{L}_\delta}(x)$, the solution of Eq.~\eqref{eq:nonlocalgreen}, has the following exact expression 
	\begin{align} 
		G_{x_0}^{\mathcal{L}_\delta}(x) = \frac{1}{\left(1-(\delta k)^2\right)^2} G_{x_0}(x) + \frac{\delta^2}{1-(\delta k)^2} \mathcal{D}(x-x_0),\label{GF}
	\end{align}
	where $G_{x_0}(x)$ is given by \eqref{eq:green} with $\tilde k =k\sqrt{1/(1-(\delta k)^2)}$.
\end{corollary}
\begin{proof} 
By definition of \eqref{eq:barqdef} and a simple calculation, at point $x_0$ we get
\begin{align}
u^w(x_0) = & \frac{1}{2\delta(1+\tilde k^2\delta^2)} \int_\R u(t+x_0) e^\frac{-|t|}{\delta} \mathrm{d}t \nonumber\\
=& \frac{1-(\delta k)^2}{2\delta} \int_\R u(t+x_0) e^\frac{-|t|}{\delta} \mathrm{d}t, 
 \label{eq:barUexpKernel}
\end{align}
where $u(x)$ is the solution to Eq. \eqref{eq:nonlocalHelmholtz}, and in the last identity we use the fact $1/(1+(\delta\tilde k)^2) = 1-(\delta k)^2$, which can be directly calculated from the relation $\tilde k = k\sqrt{1/(1-(\delta k)^2)}$. On the other hand, $u(x)$ at point $x_0$ also satisfies
\begin{align}
	\frac{1}{2\delta^3} \int_\R u(t+x_0) e^\frac{-|t|}{\delta} \mathrm{d}t &= \frac{1}{2\delta^3} \int_\R u(x_0) e^\frac{-|t|}{\delta} \mathrm{d}t - k^2 u(x_0) -f(x_0) \nonumber\\
	&= \frac{1}{\delta^2}\left[(1 - \delta^2k^2)u(x_0) - \delta^2 f(x_0) \right].
	\label{eq:barUexpKernel2}
\end{align}
From \eqref{eq:barUexpKernel2}, we obtain
\begin{align}
	u(x_0) = & \frac{1}{2\delta\left(1-(\delta k)^2\right)} \int_\R u(t+x_0) e^\frac{-|t|}{\delta} \mathrm{d}t + \frac{\delta^2}{1-(\delta k)^2} f(x_0) \notag \\
	= & \frac{1}{\left(1-(\delta k)^2\right)^2} u^w(x_0) + \frac{\delta^2}{1-(\delta k)^2} f(x_0) \nonumber\\
	= & \frac{1}{\left(1-(\delta k)^2\right)^2} \int_\R G_{x_0}(x)f(x)\mathrm{d}x + \frac{\delta^2}{1-(\delta k)^2} f(x_0) \nonumber\\
	=& \int_\R \Big(\frac{1}{\left(1-(\delta k)^2\right)^2} G_{x_0}(x)+ \frac{\delta^2}{1-(\delta k)^2}\mathcal{D}(x-x_0)\Big)f(x)\mathrm{d}x, \label{eq:uDeltaExp}
\end{align}
where the second identity uses \eqref{eq:barUexpKernel}, the third identity uses \eqref{eq:baruExp}, and the last identity uses the definition of Dirac delta function. Finally from the identity \eqref{eq:uDeltaExp}, we have that the Green's function $G_{x_0}^{\mathcal{L}_\delta}(x) $ satisfies \eqref{GF}. The proof is completed. 
\end{proof}
\begin{remark}
(a) Since the far field boundary condition for the nonlocal Helmholtz equation is still an open problem, we verify that Green's function $G_{x_0}^{\mathcal{L}_\delta}(x)$  is unique only by the property, and when the nonlocal interaction vanishes, i.e., $\delta\to0$, the Green's function $G_{x_0}^{\mathcal{L}_\delta}(x)$ converges to the Green's function $g_{x_0}(x)$ of the local Helmholtz equation with Sommerfeld radiation condition
\begin{align*}
	-\partial_x^2 g_{x_0}(x) - k^2 g_{x_0}(x) & =  D(x-x_0)     \quad & & x\in\mathbb{R},        \\
	\Big|\partial_r g_{x_0}(x) -\mathbf{i} k g_{x_0}(x)\Big| &\to  0  & & \mathrm{as}\ r=|x|\rightarrow+\infty.
\end{align*}
It is compatible with the asymptotic convergence property of nonlocal problems and can be proved simply by using the fact that $\tilde k\to k$ as $\delta\to0$.

(b) We point out that the result in Corollary \ref{th:greenexpgamma} is attractive for the study of local Helmholtz equations. For a given positive \emph{wavenumber} $\tilde k$, the solution $u_{loc}(x)$ to the following Helmholtz problem
\begin{align}
	- \partial_x^2 u_{loc}(x) - \tilde k^2 u_{loc}(x) =       & f(x),\quad x\in\R,\label{eq:localtildekHelmholtz}         \\
	\partial_r u_{loc}(x) - \i\tilde k u_{loc}(x) = & o(1),\quad \mathrm{as}\ r=|x|\rightarrow+\infty,\label{eq:localtildekHelmholtzBD}
\end{align}
can be obtained by solving the nonlocal Helmholtz equation
\begin{align}
	\mathcal{L}_\delta u(x) - k^2 u(x) = f(x), \label{eq:deltanonlocalHelmholtz}
\end{align}
where $k = \tilde k/\sqrt{1+\delta^2\tilde k^2}$ obtained by \eqref{eq:ktildekdelta}. This is to say, from Eq.~\eqref{eq:uDeltaExp} we have
\begin{align}
	u_{loc}(x) = \left( 1-(\delta k)^2\right)^2 u(x) - \delta^2 \left( 1-(\delta k)^2\right) f(x). \label{eq:utildekudelta}
\end{align}
\end{remark}

In light of Corollary~\ref{th:barq}, we now present the estimates of nonlocal solutions $u(x)$ and $u(\tilde x)$.
\begin{corollary}\label{th:expnonlocalPML}
Assume that $C$ is a constant only depending on $\Omega=(-l,l)$. Then for the solution $u(x)$ of the nonlocal Helmholtz equation \eqref{eq:nonlocalHelmholtz} and its analytic continuation $\tilde{u}(x)$ by complex coordinate transform \eqref{eq:complexStreching}, we have the following estimates:
	\begin{itemize}
		\item[(1)] If $k\delta<1$ and $f\in H^1(\R)$, it holds that 
		\begin{align}
		&\|u\|_{L^2(\Omega)}	\leq \frac{1}{k(1-(\delta k)^2)}\Big(\frac{C}{\left(1-(\delta k)^2\right)^\frac12} + \delta^2k\Big) \|f\|_{L^2(\Omega)}, \label{eq:expKernelcase1L2Sta}\\
		&|u|_{H^1(\Omega)} \leq \frac{C}{\left(1-(\delta k)^2\right)^2} \|f\|_{L^2(\Omega)} + \frac{\delta^2}{1-(\delta k)^2} |f|_{H^1(\Omega)},\label{eq:expKernelcase1H1Sta}\\
			&|\tilde u(x)|\leq C \frac{1}{\left(1-(\delta k)^2\right)^2} e^{-\frac{k}{\sqrt{1-(\delta k)^2}}|\int_0^x \sigma(t)\mathrm dt|} \|f\|_{L^2(\Omega)},\quad \mathrm{for}\ |x|>l. \label{eq:expKernelcase1ExpDecaying}
		\end{align}
		\item[(2)] If $k\delta>1$ and $f\in H^1(\R)$, it holds that 
	\begin{align}
	&\norm{u}_{L^2(\Omega)} 	\leq \frac{1}{k\left((\delta k)^2-1\right)}\left(\frac{C}{\left((\delta k)^2-1\right)^\frac12} +\delta^2 k \right)\|f\|_{L^2(\Omega)},\\
	&\abs{u}_{H^1(\Omega)}\leq \frac{C}{\left(1-(\delta k)^2\right)^2} \|f\|_{L^2(\Omega)} + \frac{\delta^2}{(\delta k)^2-1} |f|_{H^1(\Omega)}, \\
	&|u(x)| \leq \frac{C}{k\left((\delta k)^2-1\right)^\frac32} e^{-\frac{k}{\sqrt{(\delta k)^2-1}}(|x|-l)} \|f\|_{L^2(\Omega)},\qquad \mathrm{for}\ |x|>l.
\end{align}
	\end{itemize}
	\end{corollary}
	\begin{proof}
	\emph{(1)} From \eqref{eq:uDeltaExp}, $u(x)$ can be expressed by 
	\begin{align}
		u(x) = \frac{1}{\left(1-(\delta k)^2\right)^2} u^w(x) + \frac{\delta^2}{1-(\delta k)^2} f(x). \label{eq:expnonlocalPMLcase1eq1}
	\end{align}	
	 Since $\tilde k=k\sqrt{1/(1-(\delta k)^2)}$, we have $\tilde k\in \R^+$ if $\delta k<1$. Combining \eqref{eq:expnonlocalPMLcase1eq1} with Corollary~\ref{th:barq} which shows $\|u^w\|_{H^1(\Omega)} + \tilde k\|u^w\|_{L^2(\Omega)} \leq C \|f\|_{L^2(\Omega)}$, we get the stability estimates~\eqref{eq:expKernelcase1L2Sta} and \eqref{eq:expKernelcase1L2Sta}.
		
		On the other hand, since $f(x)$ is supported in $\Omega=(-l,l)$, the nonlocal PML solution $\tilde{u}(x)$ can be expressed by 
\begin{align} \label{Ux}
	\tilde u(x) = \frac{1}{\left(1-(\delta k)^2\right)^2} \tilde u^w(x)\quad \mathrm{for}\ |x|>l.
\end{align}
By combining Eq. \eqref{Ux} with \eqref{eq:nonlocalAvgExpDecayCase1} in Corollary~\ref{th:barq}, we get the estimate \eqref{eq:expKernelcase1ExpDecaying}.

	\emph{(2)} if $\delta k>1$, $\tilde k= \i |\tilde k|=\i{k}/{\sqrt{(\delta k)^2-1}}$. By Corollary~\ref{th:barq} and \eqref{eq:expnonlocalPMLcase1eq1} we have
	\begin{align*}
	\norm{u}_{L^2(\Omega)} \leq & \frac{1}{\left(1-(\delta k)^2\right)^2} \norm{u^w}_{L^2(\Omega)} + \frac{\delta^2}{(\delta k)^2-1} \|f\|_{L^2(\Omega)} \\
			\leq&\frac{1}{\left(1-(\delta k)^2\right)^2} \cdot \frac{C}{|\tilde k|} \|f\|_{L^2(\Omega)} + \frac{\delta^2}{(\delta k)^2-1} \|f\|_{L^2(\Omega)}\\
			\leq&\left(\frac{C}{k\left((\delta k)^2-1\right)^\frac32} + \frac{\delta^2}{(\delta k)^2-1} \right)\|f\|_{L^2(\Omega)}
	\end{align*}
	and
	\begin{align*}
	\abs{u}_{H^1(\Omega)} \leq & \frac{1}{\left(1-(\delta k)^2\right)^2} \abs{u^w}_{H^1(\Omega)} + \frac{\delta^2}{(\delta k)^2-1} |f|_{H^1(\Omega)} \\
			\leq&\frac{C}{\left(1-(\delta k)^2\right)^2} \|f\|_{L^2(\Omega)} + \frac{\delta^2}{(\delta k)^2-1} |f|_{H^1(\Omega)}.
		\end{align*}
Again, noting that $f(x)$ is supported in $\Omega = (-l,l)$, for $|x|\geq l$ the nonlocal Helmholtz solution $u$ satisfies 
		\begin{align}
	u(x) = \frac{1}{\left(1-(\delta k)^2\right)^2} u^w(x),
\end{align}
which combining with Corollary~\ref{th:barq} implies
\begin{align*}
	|u(x)| \leq& \frac{1}{\left(1-(\delta k)^2\right)^2} \cdot \frac{C}{|\tilde k|} e^{-\Im(\tilde k)(|x|-l)} \|f\|_{L^2(\Omega)} \\
	\leq& \frac{C}{k\left((\delta k)^2-1\right)^\frac32} e^{-\frac{k}{\sqrt{(\delta k)^2-1}}(|x|-l)} \|f\|_{L^2(\Omega)}.
\end{align*}
The proof is completed.
	\end{proof}
	\begin{figure}[htbp]
\centering
	\includegraphics[width=0.33\textwidth]{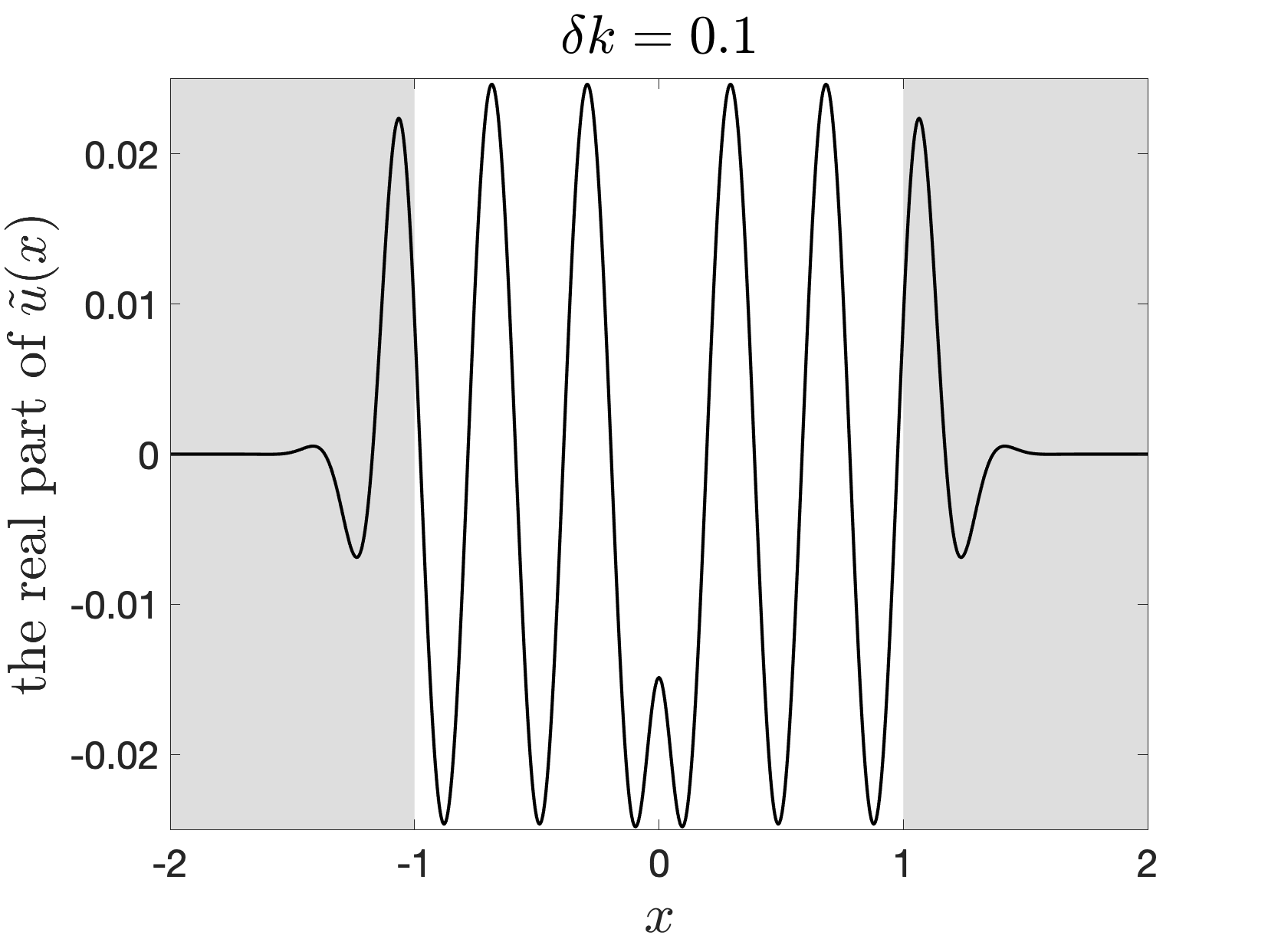}
	\includegraphics[width=0.33\textwidth]{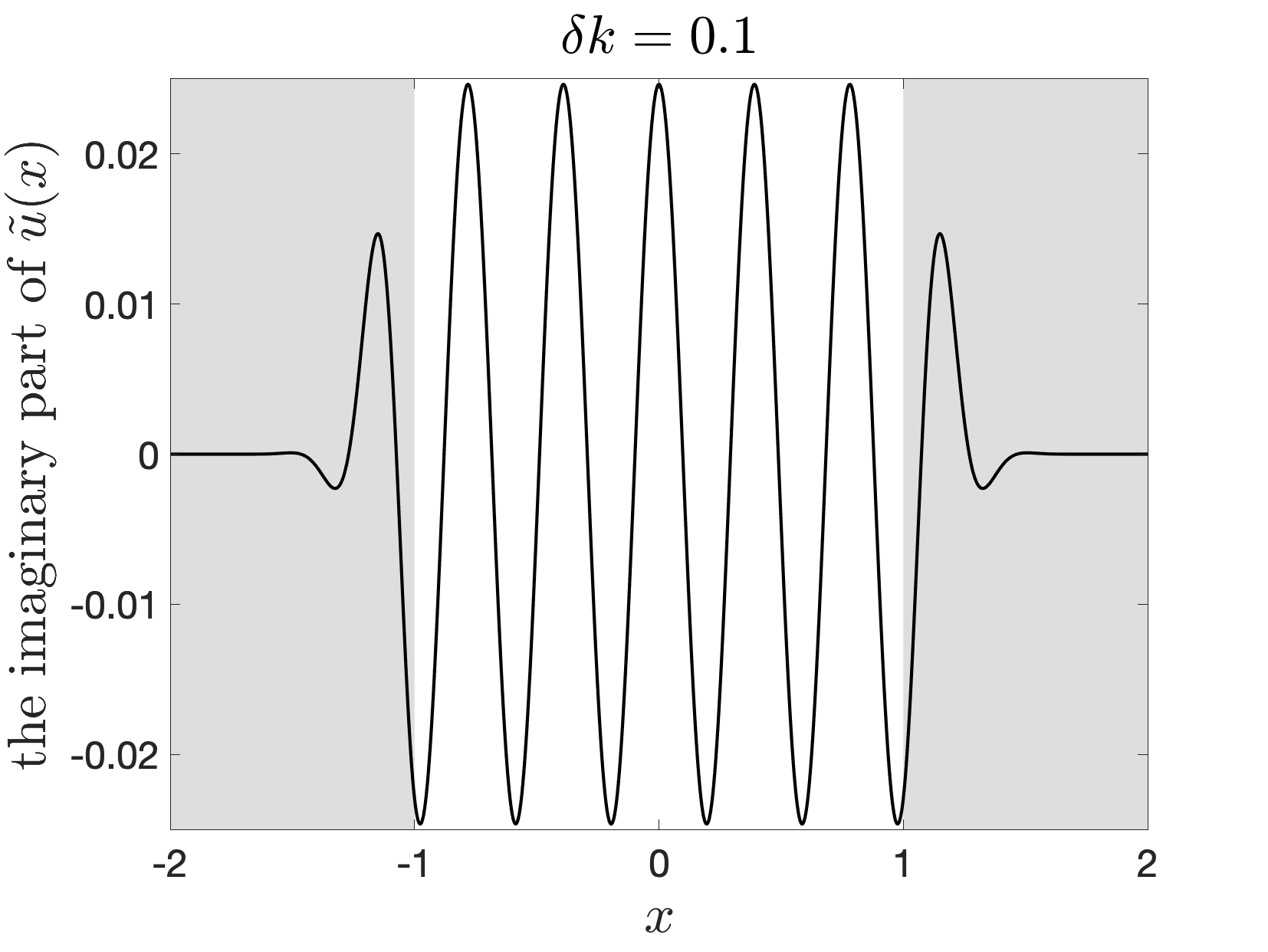}\\
	\includegraphics[width=0.33\textwidth]{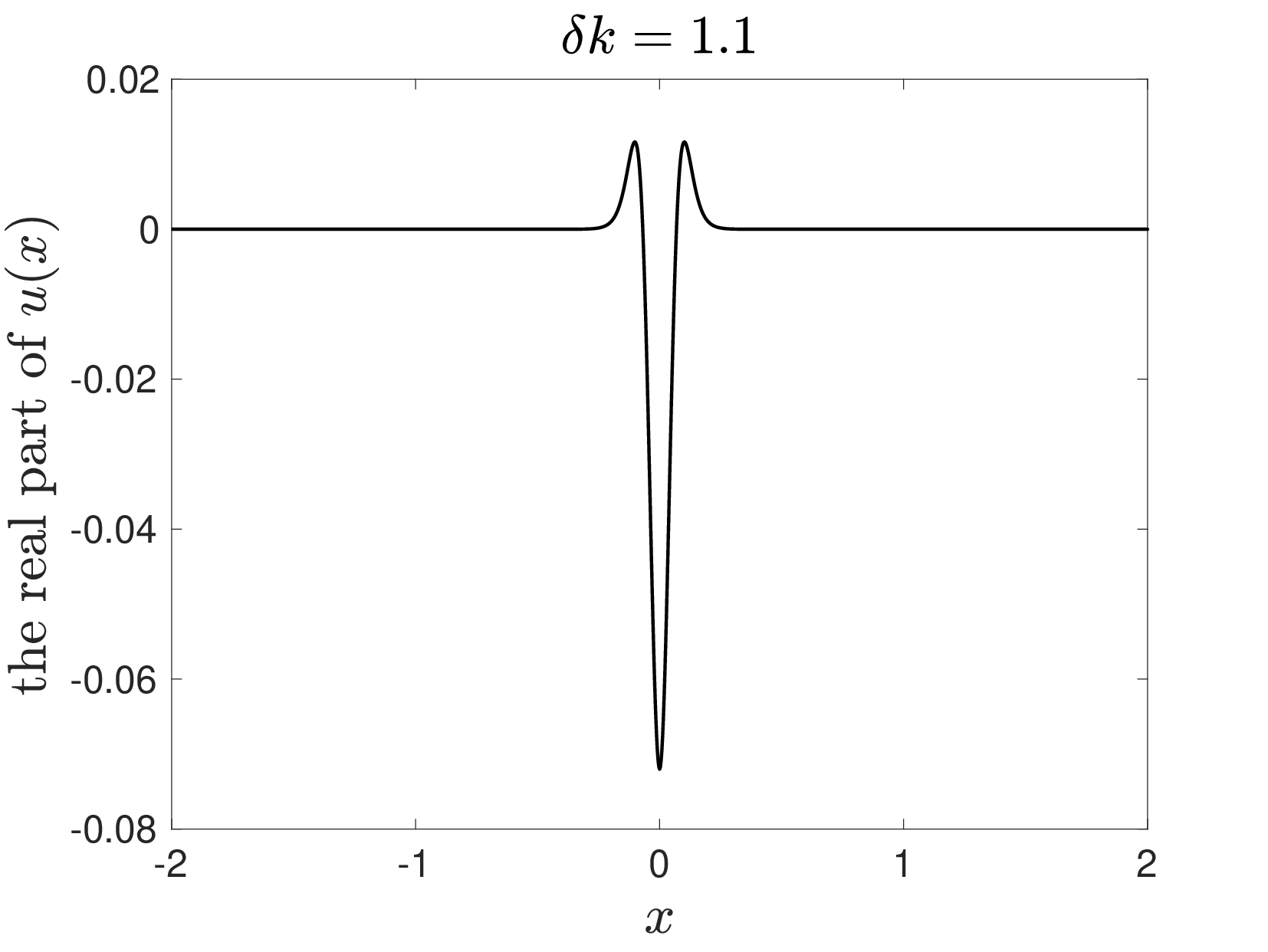}
	\includegraphics[width=0.33\textwidth]{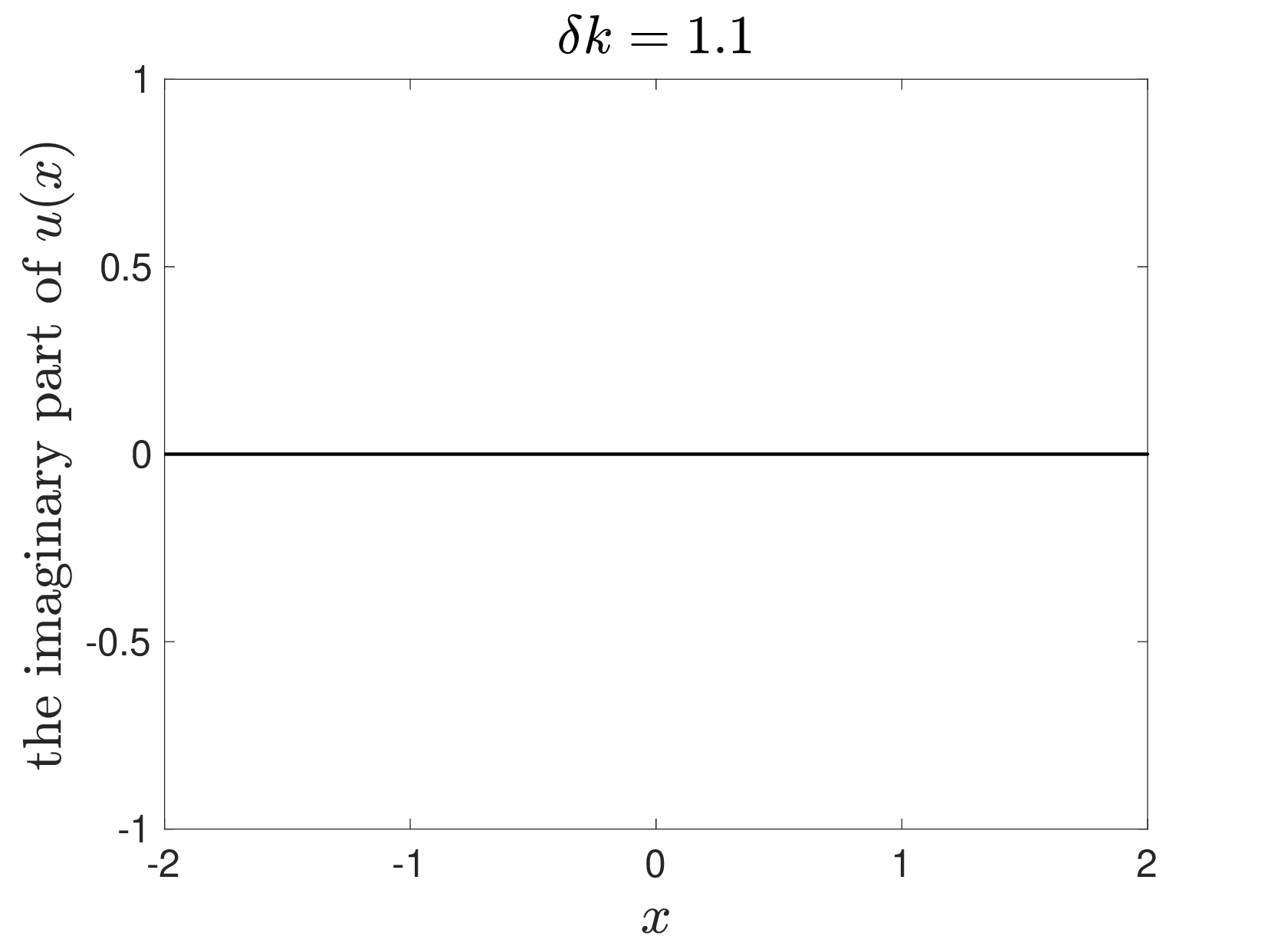}
	\caption{As an illustration of the conclusion in Corollary \ref{th:expnonlocalPML}, Top from left to right plot the real part and the imaginary part of an exact solution $\tilde u(x)$ under condition $\delta k=0.1<1$ with $k =16$. Bottom from left to right plot the real part and the imaginary part of an exact solution $u(x)$ under condition $\delta k=1.1>1$ with $k=16$. The PML layers are shaded in light grey.} \label{fig:ex1_solution}
\end{figure}	

As an illustration of the conclusion in Corollary \ref{th:expnonlocalPML}, Figure~\ref{fig:ex1_solution} shows two solutions with taking $k=16$: one is the nonlocal PML solution $\tilde{u}(x)$ under the condition that $k\delta=0.1$ and the other is the nonlocal Helmholtz solution $u(x)$ under the condition $k\delta=1.1$. One can see that $\tilde{u}(x)$ decays exponentially as $x$ goes to $\infty$ under condition $k\delta<1$, and $u(x)$ without any modification decays exponentially itself under condition $k\delta>1$.

\subsection{The truncated nonlocal problems on bounded domains}
Based on the analysis in corollaries \ref{th:barq} and \ref{th:expnonlocalPML}, we can conclude two folds for a fixed $\delta$: (i) for \emph{wavenumber} $k$ such that $\tilde k\in \R^+$, 
the nonlocal PML solution decays exponentially in the PML layer; (ii) for \emph{wavenumber} $k$ large enough such that $\tilde k\in\mathbb{C}$ with $\Im(\tilde k)>0$, the nonlocal Helmholtz solution itself decays exponentially outside the computational domain of interest. Hence the results in corollaries \ref{th:barq} and \ref{th:expnonlocalPML} suggest that we can consider two cases to \emph{truncate} the computational region at some sufficiently large $x$ by putting forced boundary conditions, such as homogeneous Dirichlet boundary conditions.

\begin{itemize}
	\item[] {\bf Case 1.} For  \emph{wavenumber} $k$ such that $\tilde k\in \R^+$, we need to truncate the nonlocal PML equation \eqref{eq:PMLeq}, namely, let the truncated PML solution $\hat{\tilde u}(x)$ satisfy 
	  \begin{align}
		  \tilde{\mathcal{L}}_\delta\hat{\tilde u}(x) - k^2\omega(x)\hat{\tilde u}(x) = & f(x),\quad & & |x|\leq l+d,   \label{eq:truPML1}   \\
		  \hat{\tilde u}(x) =             & 0,   & & l+d< |x|\leq l+d+l_\gamma,
		  \label{eq:truPML2}
	  \end{align}
	  where $d$ represents the thickness of PML layers.

	\item[] {\bf Case 2.} For  \emph{wavenumber} $k$ large enough such that $\tilde k\in\mathbb{C}$ with $\Im(\tilde k)>0$, we only need to truncate the nonlocal Helmholtz equation~\eqref{eq:nonlocalHelmholtz} by imposing a homogeneous Dirichlet boundary condition surrounding $\Omega$, namely, the truncated Helmholtz solution $\hat u(x)$ is governed by 
	  \begin{align}
		  \mathcal{L}_\delta\hat u(x) - k^2\omega(x)\hat u(x) = & f(x),\quad & &|x|\leq l+d,  \label{eq:trueq1}  \\
		  \hat u(x) =         & 0,   & & l+d< |x|\leq l+d+l_\gamma,
		  \label{eq:trueq2}
	  \end{align}
	  where $d$ is the thickness of layers on which the solution decays exponentially.
\end{itemize}

\begin{remark}

In the previous section, we show that for the kernel $\gamma_1( s)=\frac{1}{2}e^{-| s|}$, the solution to the local Helmholtz problem~\eqref{eq:localtildekHelmholtz}--\eqref{eq:localtildekHelmholtzBD} can be obtained by solving the nonlocal Helmholtz equation~\ref{eq:deltanonlocalHelmholtz} (see \eqref{eq:utildekudelta}). Similarly, for positive $\tilde k$ the solution of the following truncated local PML problem
\begin{align*}
	-\partial_x \left( \frac{1}{\omega(x)}\partial_x\hat{\tilde u}_{\tilde k}(x)\right) - \tilde k^2 \omega(x) \hat{\tilde u}_{\tilde k}(x) = & f(x),\quad & & x\in (-l-d,l+d),   \\
	\hat{\tilde u}_{\tilde k}(x)=                     & 0,\quad & & x=-l-d\ \mathrm{or}\ l+d,
\end{align*}
can be approximated by using the solution $\hat{\tilde u}$ of \eqref{eq:truPML1}--\eqref{eq:truPML2} with the kernel $\gamma_1=\frac12 e^{-| s|}$, due to the property that all of these solutions in the whole space decay exponentailly as $|x|\rightarrow+\infty$.
\end{remark}

\section{Discretization}\label{sec:Discretization}
We have reformulated the model problem \eqref{eq:nonlocalHelmholtz} into two kinds of truncated problems, i.e., the problem \eqref{eq:truPML1}-\eqref{eq:truPML2} and problem \eqref{eq:trueq1}-\eqref{eq:trueq2}, and here present the corresponding discrete schemes. As discussed in \cite{du2012analysis, zhou2010mathematical, du2019nonlocal, du2018nonlocal}, the nonlocal operator converges to the corresponding local operator as the nonlocal interaction horizon $\delta$ vanishes. Hence, it is useful to use the asymptotic compatibility (AC) scheme, a concept developed in \cite{TianDu}, to discretize the nonlocal operator as AC schemes can preserve the analogous limit \eqref{eq:NTL} in a discrete level as both mesh size $h$ and horizon $\delta \rightarrow 0$. Here we mainly use the AC scheme given in \cite{du2018nonlocal} to discretize the nonlocal operator. To do so, we assume that the domain $[-l-d-l_\gamma,l+d+l_\gamma]$ has been discretized by a uniform grid $x_{-N}, x_{-N+1},\cdots,x_N$ with spacing $h$, and there exists an integer $M$ such that $l+d=Mh$ for simplicity. Under the framework of the AC scheme in \cite{du2018nonlocal}, we take $\gamma(x,y)=\gamma_\delta(x-y)$ and $\tilde \gamma(x,y)=\gamma_\delta(\tilde x,\tilde y)\omega(x)\omega(y)$, and define
\begin{align*}
F(x,y,s) := \frac{\hat u(x)-\hat u(y)}{x-y} s \gamma\left( \frac{x+y}{2}-\frac{s}{2},\frac{x+y}{2}+\frac{s}{2} \right),\\
\tilde F(x,y,s) := \frac{\hat{\tilde u}(x)-\hat{\tilde u}(y)}{y-x} s \tilde \gamma\left(\frac{x+y}{2}-\frac{s}{2},\frac{x+y}{2}+\frac{s}{2}\right).
\end{align*}
It's clear that
\begin{align*}
\mathcal{L}_\delta\hat u(x) := \int_\R F(x,y,y-x) \mathrm{d}y,\\
\tilde{\mathcal{L}}_\delta\hat{\tilde u}(x) := \int_{\R} \tilde F(x,y,y-x) \mathrm{d}y.
\end{align*}
Then we expand $F(x,y,s)$ and $\tilde F(x,y,s)$ in test functions with respect to $y$
\begin{align*}
F_h(x,y,s) =& \sum_m \phi_m(y) F(x,x_m,s),\\
\tilde F_h(x,y,s) =& \sum_m \phi_m(y) \tilde F(x,x_m,s),
\end{align*}
where $\phi_m(y)$ is the hat function of width $h$ centered at $y_m=mh$. Thus,
the discrete formulation of $\mathcal{L}_\delta$ is given by
\begin{align*}
	\mathcal{L}_\delta^h\hat u(x_n) = \sum_{m\in\mathbb{Z}} a_{n,m} \hat u(x_m)
\end{align*}
with $ x_{\frac{m+n}{2}}=(x_m+x_n)/2$ and 
\begin{align*}
	a_{n,m} = \left\{
	\begin{aligned}
		&- \frac{1}{(m-n)h} \int_{\R} \Big[\phi_{m-n}(s) s \gamma\Big(x_{\frac{m+n}{2}}-\frac{s}{2},x_{\frac{m+n}{2}}+\frac{s}{2}\Big) \Big]\ds, & & m\neq n, \\
		 & -\sum_{m\neq n} a_{n,m},\quad                         & & m=n.
	\end{aligned}\right.
\end{align*}
Similarly, we give the discrete formulation of $\tilde{\mathcal{L}}_\delta$ by
\begin{align*}
	\tilde{\mathcal{L}}_\delta^h\hat{\tilde u}(x_n) = \sum_{m\in\mathbb{Z}} \tilde a_{n,m} \hat{\tilde u}(x_m),
\end{align*}
where
\begin{align*}
	\tilde a_{n,m} = \left\{
	\begin{aligned}
		&- \frac{1}{(m-n)h} \int_{\R} \Big[\phi_{m-n}(s) s \tilde \gamma\Big(x_{\frac{m+n}{2}}-\frac{s}{2},x_{\frac{m+n}{2}}+\frac{s}{2}\Big) \Big]\ds, & & m\neq n, \\
		 & -\sum_{m\neq n} \tilde a_{n,m}\quad                         & & m=n.
	\end{aligned}\right.
\end{align*}
%
Finally, the truncated problems \eqref{eq:truPML1}--\eqref{eq:truPML2} and \eqref{eq:trueq1}--\eqref{eq:trueq2} are respectively discretized by
\begin{itemize}
\item the discretization for the case with  \emph{wavenumber} $k$ such that $\tilde k\in \R^+$:
\begin{align}
\tilde{\mathcal{L}}_\delta^h \hat{\tilde u}_h(x_n) - k^2\omega(x_n)\hat{\tilde u}_h(x_n) =& f(x_n),\quad -M<n<M,\label{eq:DAScase1eq1}\\
\hat{\tilde u}_h(x_n) = & 0, \qquad \quad M\leq |n|\leq N;\label{eq:DAScase1eq2}
\end{align}
\item the discretization for the case with  \emph{wavenumber} $k$ such that $\tilde k\in\mathbb{C}$ with $\Im(\tilde k)>0$:
\begin{align}
\mathcal{L}_\delta^h \hat{u}_h(x_n) - k^2\omega(x_n)\hat{u}_h(x_n) =& f(x_n),\quad M<n<M,\label{eq:DAScase2eq1}\\
\hat{u}_h(x_n) = & 0, \qquad \quad M\leq |n|\leq N.\label{eq:DAScase2eq2}
\end{align}
\end{itemize}

%

\section{Numerical examples}\label{sec:numericalTests}
We now provide numerical examples to illustrate our mathematical analysis for the nonlocal solution behaviors, and the effectiveness of our PML strategies. In subsection~\ref{subsec:PMLmod4kernel}, we investigate the necessity to analytically continue the kernel function. We use subsection~\ref{subsec:EfficiencyPML} to investigate the efficiency of our PML technique, and the dependence of truncation errors on the width $d$ of PML layers and PML medium parameters $\sigma$. In subsection~\ref{subsec:EfficiencyScheme}, we focus on the convergence rates of numerical schemes in both $L^2$-norm and $H^1$-seminorm. In subsection~\ref{subsec:AC}, we simply illustrate the asymptotic compatibility in the discrete level.

The relative errors in both $L^2$-norm and $H^1$-seminorm over $\Omega$ are defined by
\begin{align*}
e_{L^2}(v_1,v_2)=\frac{\|v_1-v_2\|_{L^2(\Omega)}}{\|v_2\|_{H^1(\Omega)}}, \quad\quad e_{H^1}(v_1,v_2)=\frac{|v_1-v_2|_{H^1(\Omega)}}{\|v_2\|_{H^1(\Omega)}}.
\end{align*}
In the calculation, we set $l=10$ and set the PML medium as 
$$\sigma(t)=\frac{\sigma_0}{d}(|t|-l), \quad \text{for} \quad |t|>l,$$
 where $\sigma_0$ represents the PML medium parameter. The source $f$ is chosen as the Gaussian function 
\begin{align}
	f(x) = e^{-(\frac{2k}{5\pi})^2x^2}. \label{eq:narrowGaussianSource}
\end{align}
and the exponential and Gaussian kernels are taken and given by 
\begin{align}
\gamma_1^{exp}( s)=\frac12 e^{-| s|},\quad \gamma_1^{gau}( s) = \frac{4}{\sqrt{\pi}} e^{- s^2}. 
\end{align}

Noting that the kernel function in our nonlocal PML equation must be the analytic continuation of the original kernel, we here present how to implement the complex coordinate stretching for the two kernel functions above.

For $\gamma_1^{exp}(s)$ with $ s=|x-y|$ representing the distance between two location points $x$ and $y$ in the real space $\R$, we introduce the notation $\rho(x,y)=\sqrt{(x-y)^2}$. For the complex valued points in PML modifications, we define $\rho(\tilde x,\tilde y)=\sqrt{(\tilde x-\tilde y)^2}$ by taking the analytic branch such that $\Re\rho(\tilde x,\tilde y)\geq0$. Thus, for the kernels $\gamma_1^{exp}$ and $\gamma_1^{gau}$, the complex coordinate stretching is given by
\begin{align*}
	\gamma_1^{exp}(\tilde x-\tilde y) = \frac{1}{2} e^{-\rho(\tilde x,\tilde y)}\quad \text{and} \quad \gamma_1^{gau}(\tilde x-\tilde y)= \frac{4}{\sqrt{\pi}} e^{-(\tilde x-\tilde y)^2}. 
\end{align*}

In the forthcoming numerical tests, the exact nonlocal Helmholtz solution $u$ for the kernel $\gamma_1^{exp}(s)$ is obtained by using its Green's function $G_{x_0}^{\mathcal{L}_\delta}(x)$  in \eqref{GF}, and the reference solution $u$ for $\gamma_1^{gau}(s)$ is obtained by the AC schemes \eqref{eq:DAScase1eq1}--\eqref{eq:DAScase1eq2} and \eqref{eq:DAScase2eq1}--\eqref{eq:DAScase2eq2} over a domain large enough with sufficiently fine meshes.

\subsection{PML modifications for kernel functions} \label{subsec:PMLmod4kernel}
\begin{figure}
\centering
	\includegraphics[width=0.33\textwidth]{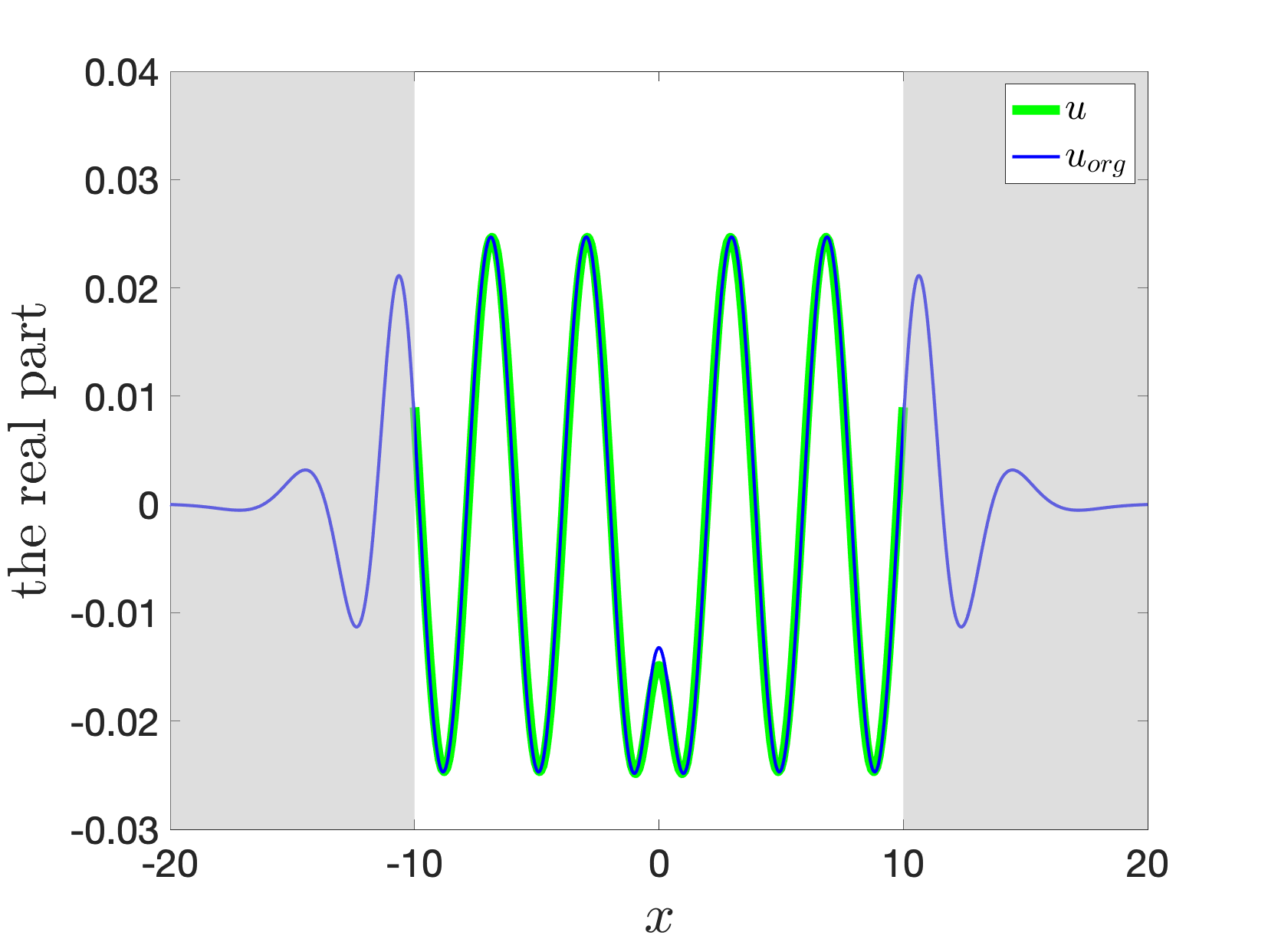}
	\includegraphics[width=0.33\textwidth]{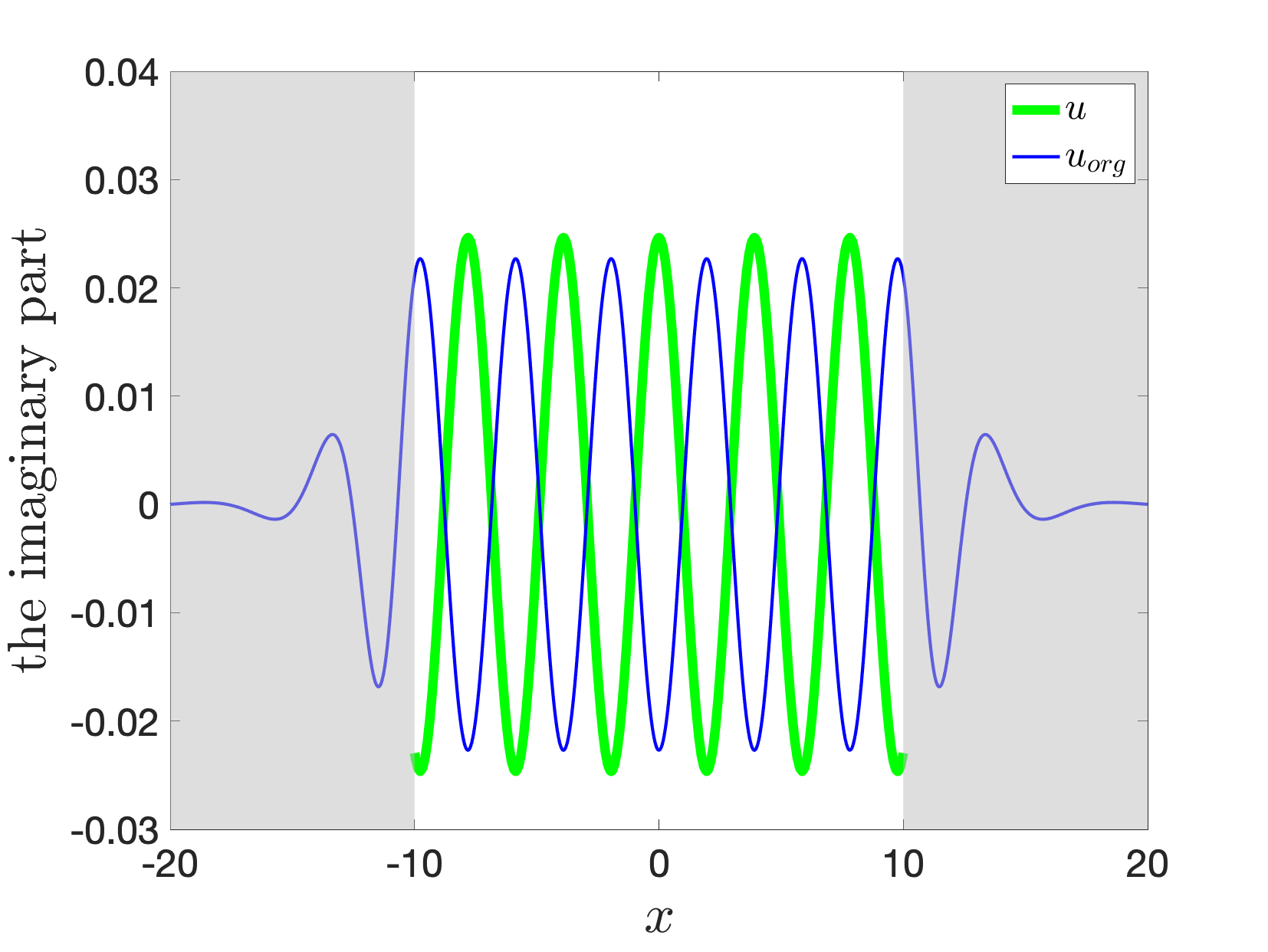}
	\caption{(Example 1) The comparison between solution $u_{org}$ of \eqref{eq:nonlocalPMLwithOriginKernel} and solution $u$ of ~\eqref{eq:nonlocalHelmholtz}. } 
	\label{fig:ex1_solution_diffkernel}
\end{figure}

{\bf Example 1.} We use this example to show that it is necessary to analytically continue the original kernel in our PML. To do so, we consider the nonlocal ``PML'' equation without modifying the kernel, namely, 
\begin{align}
\int_\R \big( u_{org}(x)-u_{org}(y) \big) \gamma_\delta(y-x) \omega(x)\omega(y) \mathrm{d}y - k^2\omega(x)u_{org}(x)=f(x). \label{eq:nonlocalPMLwithOriginKernel}
\end{align}
In the simulation, we use the kernel $\gamma_1^{exp}$, and take the parameters $k=1.6$, $\delta=1/16$. The solution $u_{org}$ of Eq.~\eqref{eq:nonlocalPMLwithOriginKernel} is obtained on a very fine mesh with $M=4096$ and $\sigma_0$ is chosen such that the exponential decaying factor (see Eq.~\eqref{eq:expKernelcase1ExpDecaying})
\begin{align*}
\frac{1}{\left(1-(\delta k)^2\right)^2} e^{-\frac{k}{\sqrt{1-(\delta k)^2}}|\int_0^{l+d} \sigma(t)\mathrm dt|}\leq 10^{-10}.
\end{align*}
Figure~\ref{fig:ex1_solution_diffkernel} shows the comparison between $u_{org}$ and the exact solution $u$ of Eq.~\eqref{eq:nonlocalHelmholtz}. One can see that the solution $u_{org}(x)$ is different from the correct solution $u$ of Eq.~\eqref{eq:nonlocalHelmholtz} on $\Omega$.

\begin{figure}
\centering
		\subfigure[]{\includegraphics[width=0.31\textwidth]{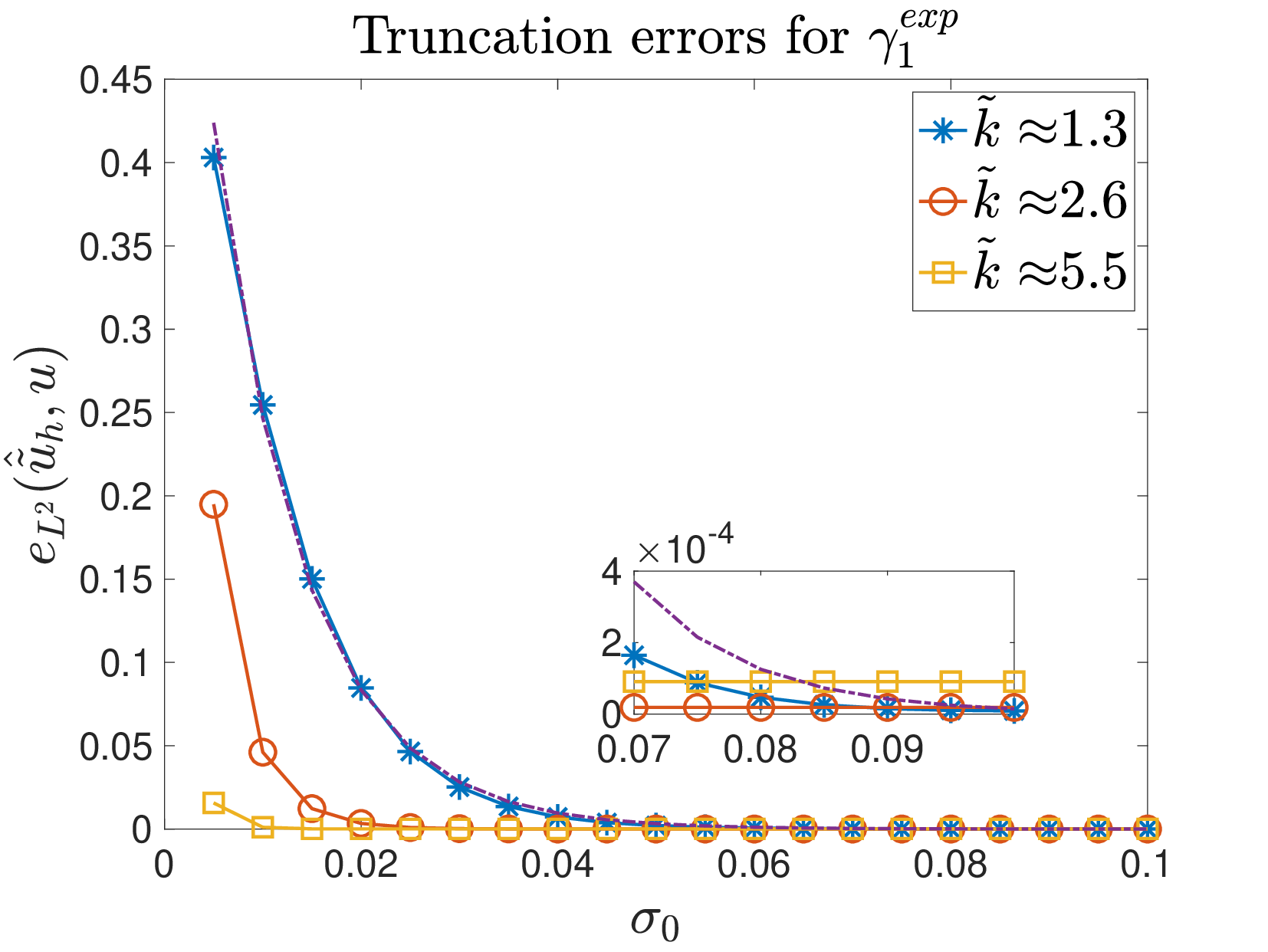}\label{subfig:ex1_truncationErrors_a}}
		\subfigure[]{\includegraphics[width=0.31\textwidth]{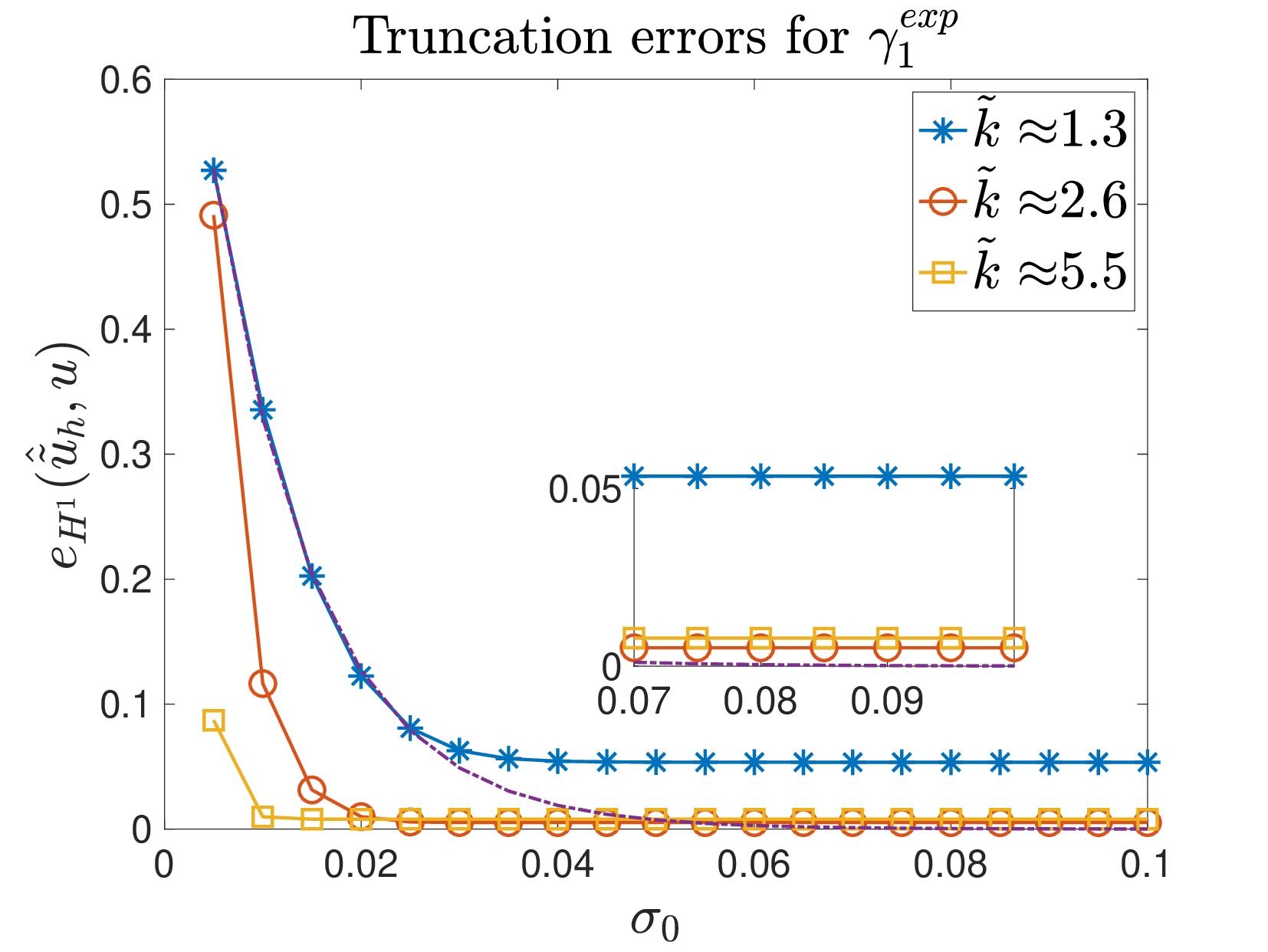}\label{subfig:ex1_truncationErrors_b}}\\
		\subfigure[]{\includegraphics[width=0.31\textwidth]{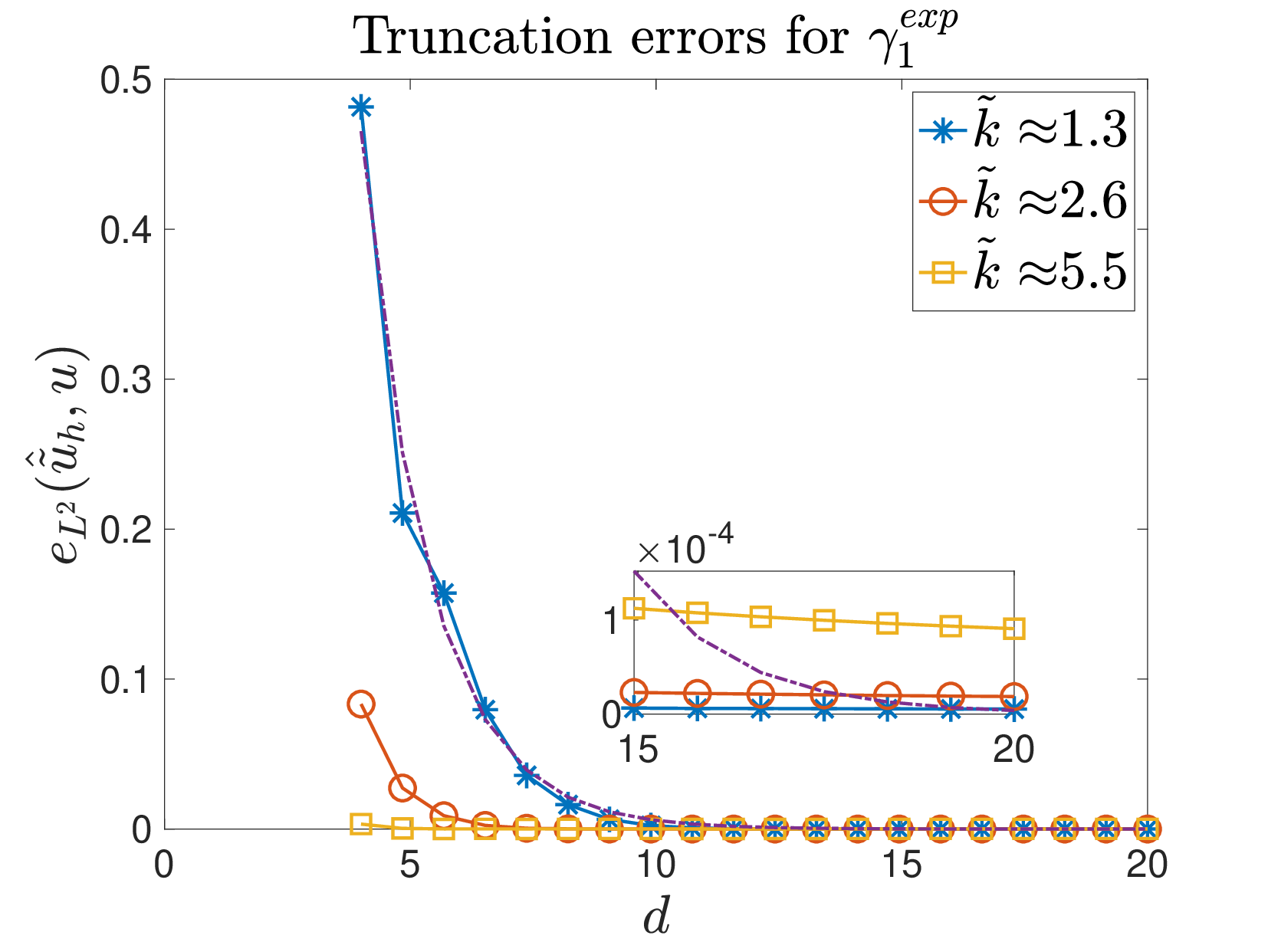}\label{subfig:ex1_truncationErrors_c}}
		\subfigure[]{\includegraphics[width=0.31\textwidth]{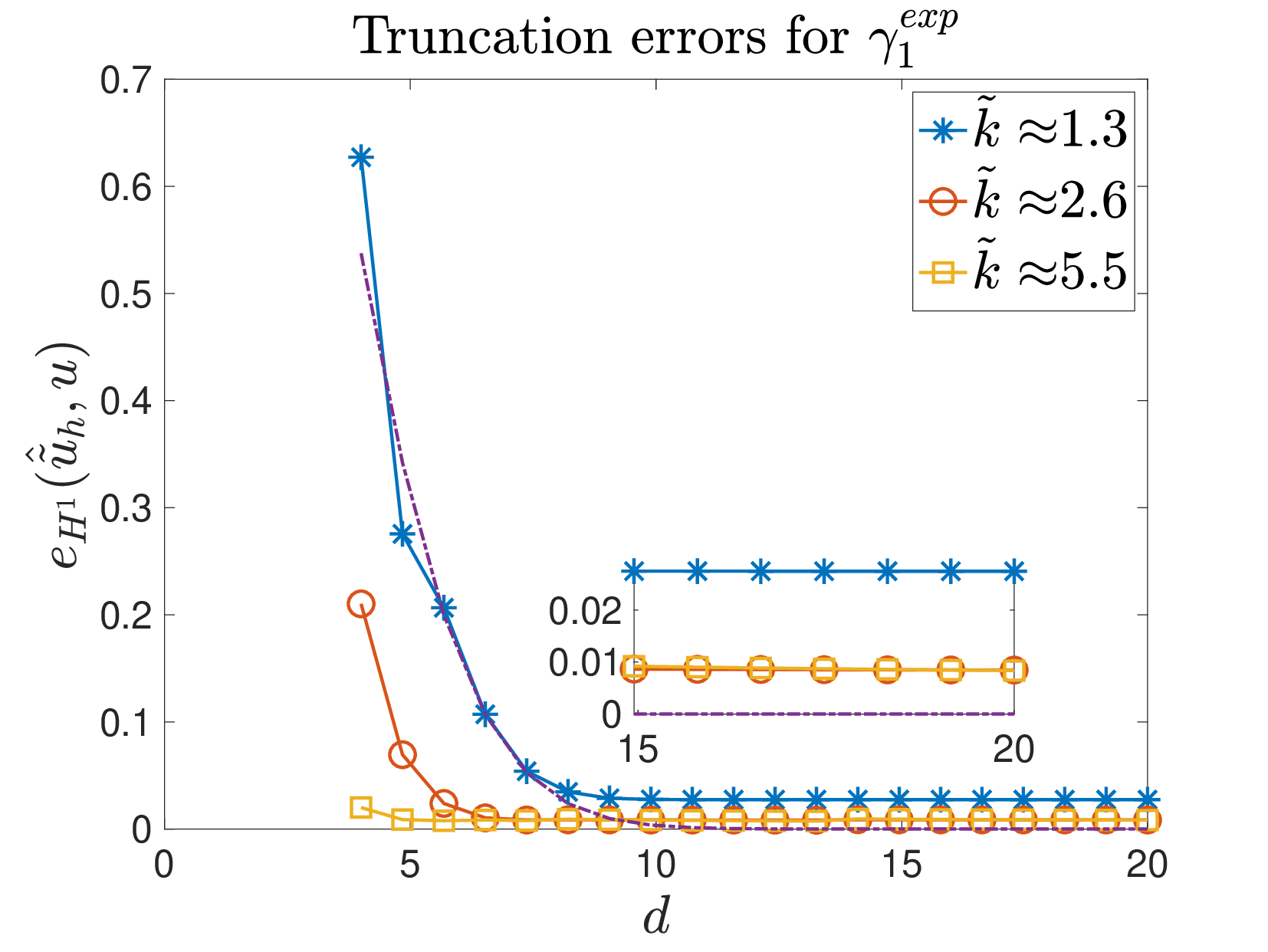}\label{subfig:ex1_truncationErrors_d}}		
	\caption{(Example 2) Truncation errors for the kernel $\gamma_1^{exp}$: (a)--(b) show errors with a fixed $d=10$ by changing $\sigma_0$, and (c)--(d) show errors with a fixed $\sigma_0=0.05$ by changing $d$. Dash-dotted lines are reference lines plotted by using $y(\sigma_0)=c_1e^{-c_2\sigma_0}$ and $y(d)=c_1e^{-c_2d}$ for some constants $c_1$ and $c_2$, respectively.
	} \label{fig:ex1_truncationErrors}
\end{figure}

\subsection{Truncation errors of the nonlocal PML solutions} \label{subsec:EfficiencyPML}
{\bf Example 2.} Here we show that the truncation errors between the nonlocal Helmholtz solution $u$ and the truncated PML solution $\hat{\tilde u}$ decay exponentially by investigating the dependence of these truncation errors on the width $d$ of PML layers and the PML parameter $\sigma_0$. In the simulations, we take the \emph{wavenumbers} $k=2\pi/5$, $4\pi/5, 8\pi/5$ and $\delta=\frac{1}{4\pi}$, which implies $\tilde k\approx1.3, 2.5$ and $5.2$ by Eq.~\eqref{eq:gausskerneltildek}. In practical computations, the truncated PML solution $\hat{\tilde u}$ is replaced by its discrete solution $\hat{\tilde u}_h$ with $h=40/8056$.

For the kernel $\gamma_1^{exp}$, Figure~\ref{fig:ex1_truncationErrors} plots the truncation errors by changing $\sigma_0$ and $d$, respectively. Specially, we fix $d=10$ in subfigures~\ref{subfig:ex1_truncationErrors_a}--\ref{subfig:ex1_truncationErrors_b} by taking different values of $\sigma_0$, and fix $\sigma_0=0.05$ in subfigures~\ref{subfig:ex1_truncationErrors_c}--\ref{subfig:ex1_truncationErrors_d} by taking different values of $d$. 

\begin{figure}
\centering
		\subfigure[]{\includegraphics[width=0.31\textwidth]{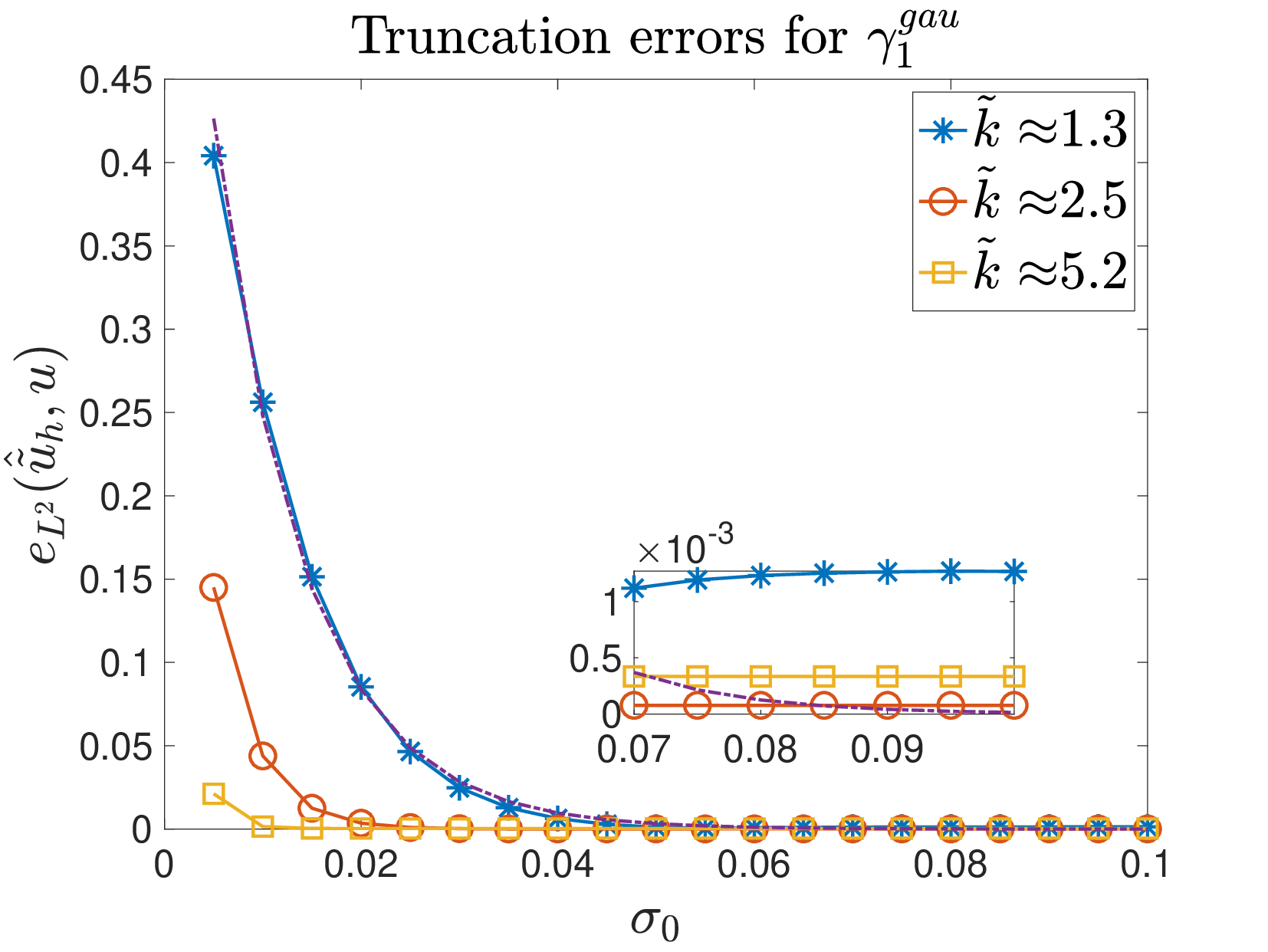}\label{subfig:ex2_truncationErrors_a}}
		\subfigure[]{\includegraphics[width=0.31\textwidth]{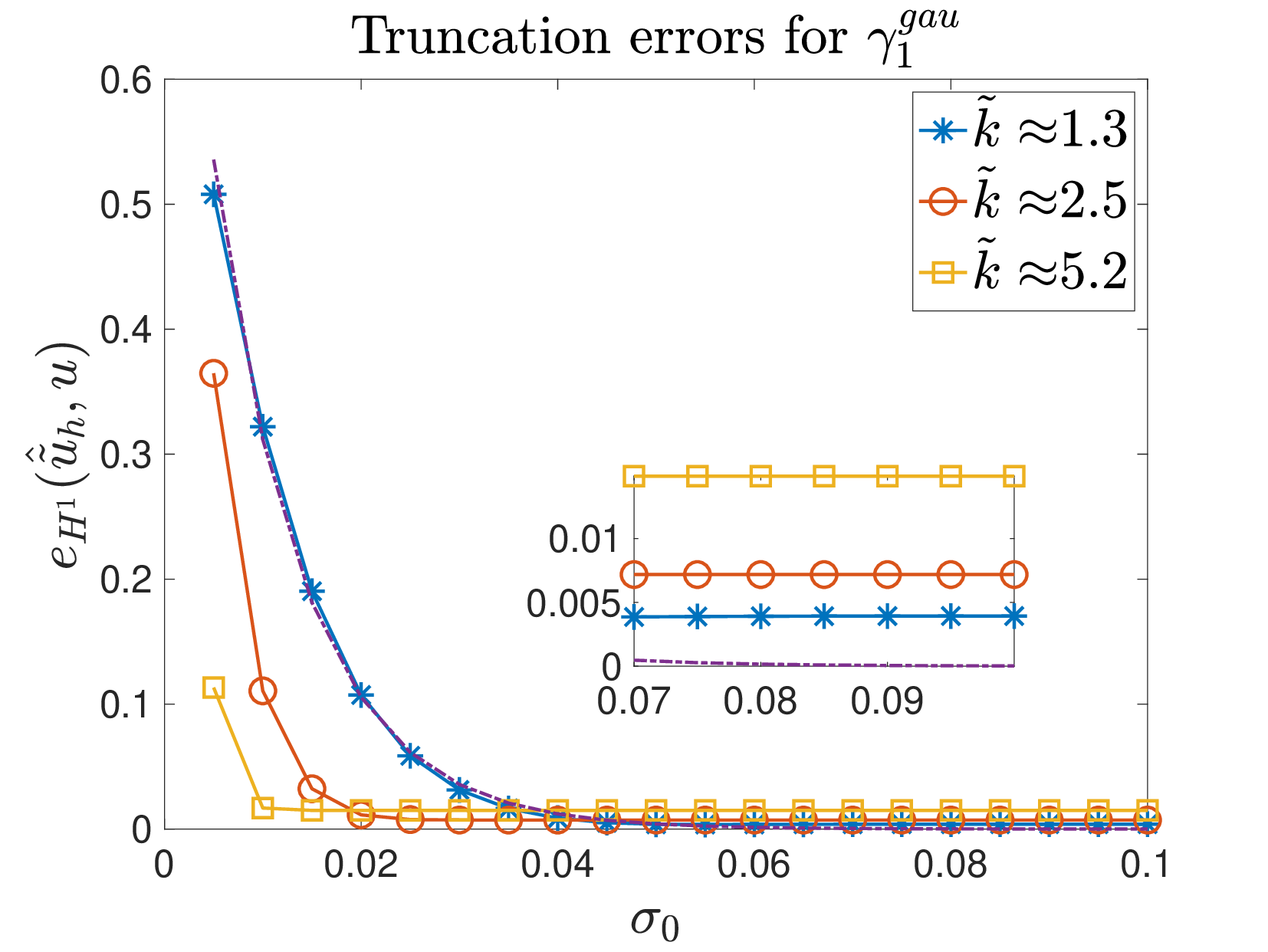}\label{subfig:ex2_truncationErrors_b}}\\
		\subfigure[]{\includegraphics[width=0.31\textwidth]{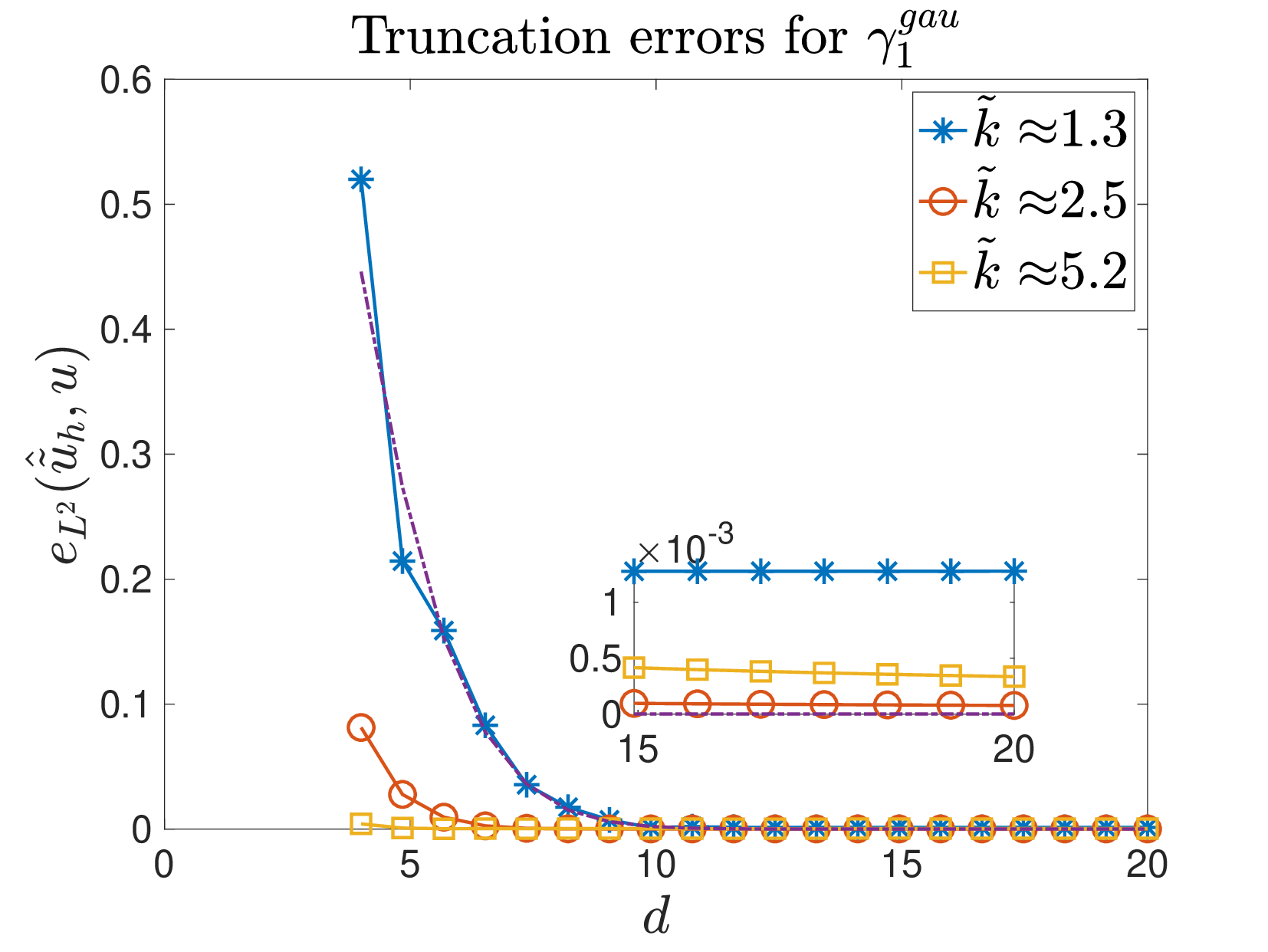}\label{subfig:ex2_truncationErrors_c}}
		\subfigure[]{\includegraphics[width=0.31\textwidth]{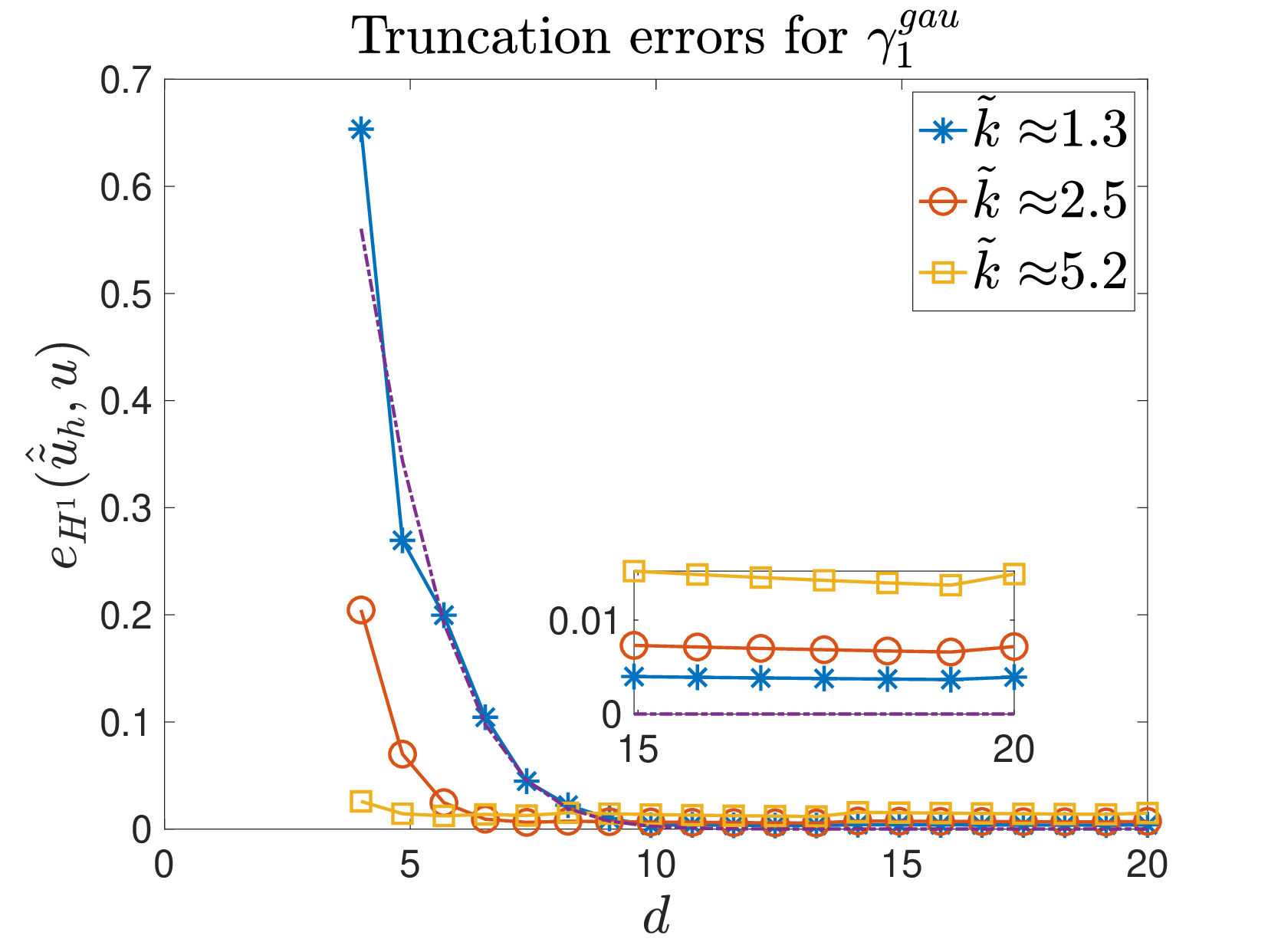}\label{subfig:ex2_truncationErrors_d}}
	\caption{(Example 2) Truncation errors for the kernel $\gamma_1^{gau}$: (a)--(b) show errors with a fixed $d=10$ by changing $\sigma_0$, and (c)--(d) show errors with a fixed $\sigma_0=0.05$ by changing $d$. Dash-dotted lines are reference lines plotted by using $y(\sigma_0)=c_1e^{-c_2\sigma_0}$ and $y(d)=c_1e^{-c_2d}$ for some constants $c_1$ and $c_2$, respectively. } \label{fig:ex2_truncationErrors}
\end{figure}

Similarly, for $\gamma_1^{gau}$, Figure~\ref{fig:ex2_truncationErrors} plots the truncation errors by changing $\sigma_0$ and $d$, respectively. We fix $d=10$ in subfigures~\ref{subfig:ex2_truncationErrors_a}--\ref{subfig:ex2_truncationErrors_b}, and fix $\sigma_0=0.05$ in subfigures~\ref{subfig:ex2_truncationErrors_c}--\ref{subfig:ex2_truncationErrors_d}. 

From Figures \ref{fig:ex1_truncationErrors} and \ref{fig:ex2_truncationErrors}, one can observe that the truncation errors decay exponentially. We point out that the limits of these errors are not zero, while they should be zero by our theoretical findings. The reason is that we replace the truncated PML solutions by their corresponding discrete solutions in the simulations, hence the limits are actually the discrete errors of the AC scheme. This is to say, the limits will decease when using smaller mesh size $h$.

\subsection{Convergence rates of the discrete schemes} \label{subsec:EfficiencyScheme}
{\bf Example 3.} We here investigate the convergence rates of discrete solutions of schemes~\eqref{eq:DAScase1eq1}--\eqref{eq:DAScase1eq2} and \eqref{eq:DAScase2eq1}--\eqref{eq:DAScase2eq2} in the relative $L^2$-norm and $H^1$-seminorm.

In the simulations, we take $d=10$ and choose $\sigma_0$ such that the decay factor (see \eqref{eq:nonlocalAvgExpDecayCase1}) satisfies
\begin{align*}
 e^{-\tilde k|\int_0^{l+d} \sigma(t)\mathrm dt|} \leq 10^{-10}.
\end{align*}

\begin{figure} 
	\centering
		\subfigure[]{\includegraphics[width=0.33\textwidth]{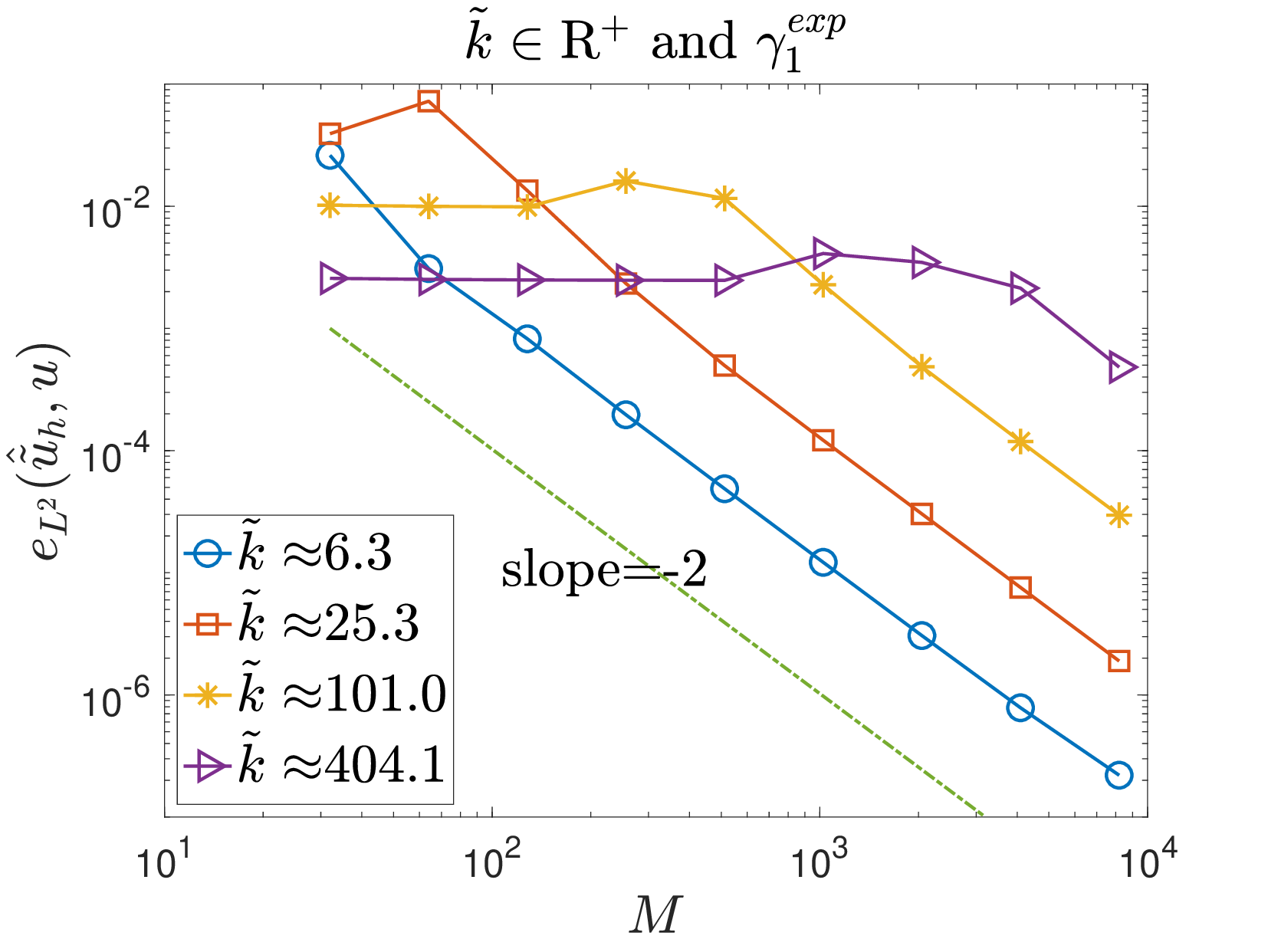}\label{subfig:ex1_fixed_deltaXk_a}}
		\subfigure[]{\includegraphics[width=0.33\textwidth]{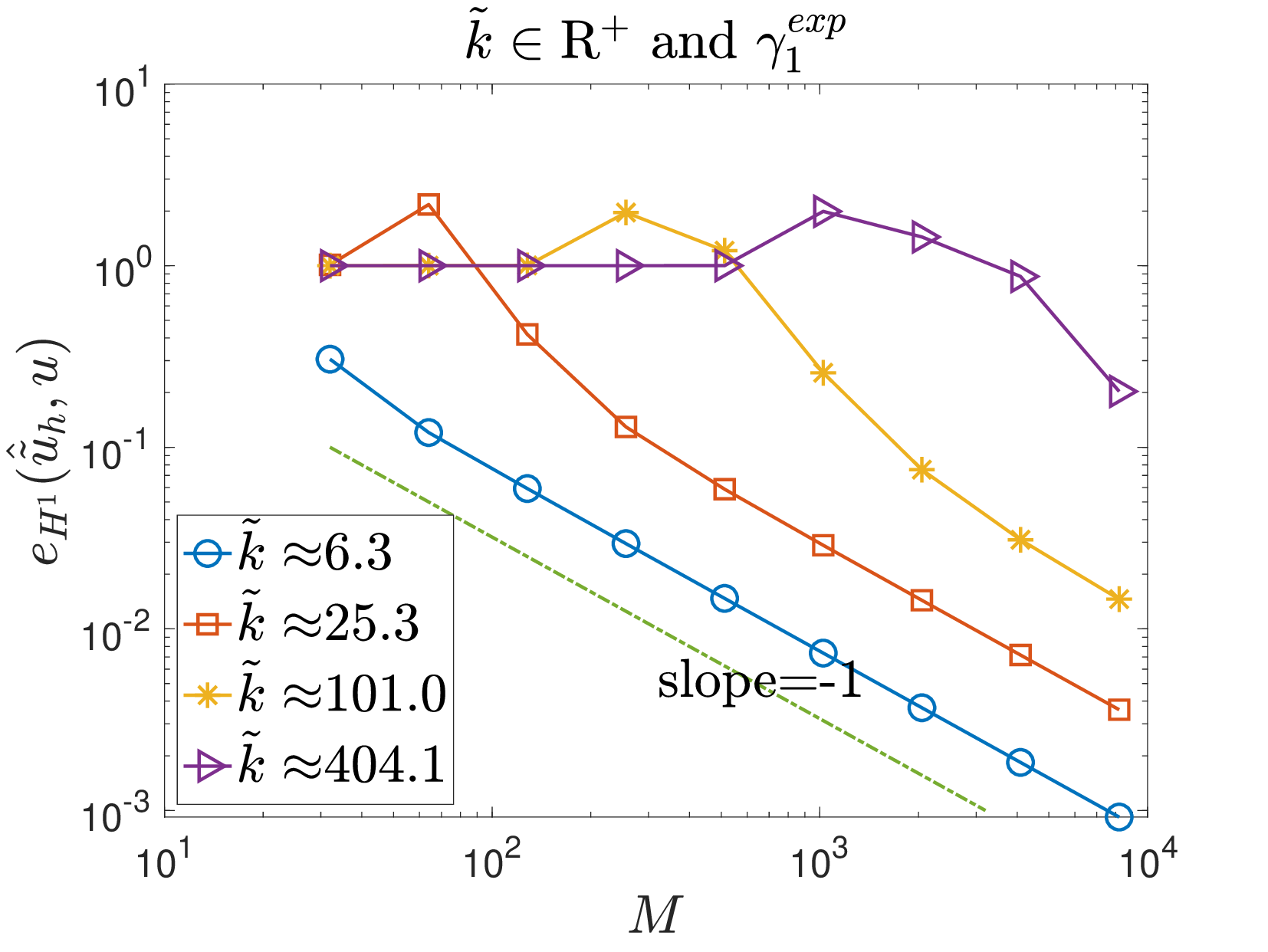}\label{subfig:ex1_fixed_deltaXk_b}}\\
		\subfigure[]{\includegraphics[width=0.33\textwidth]{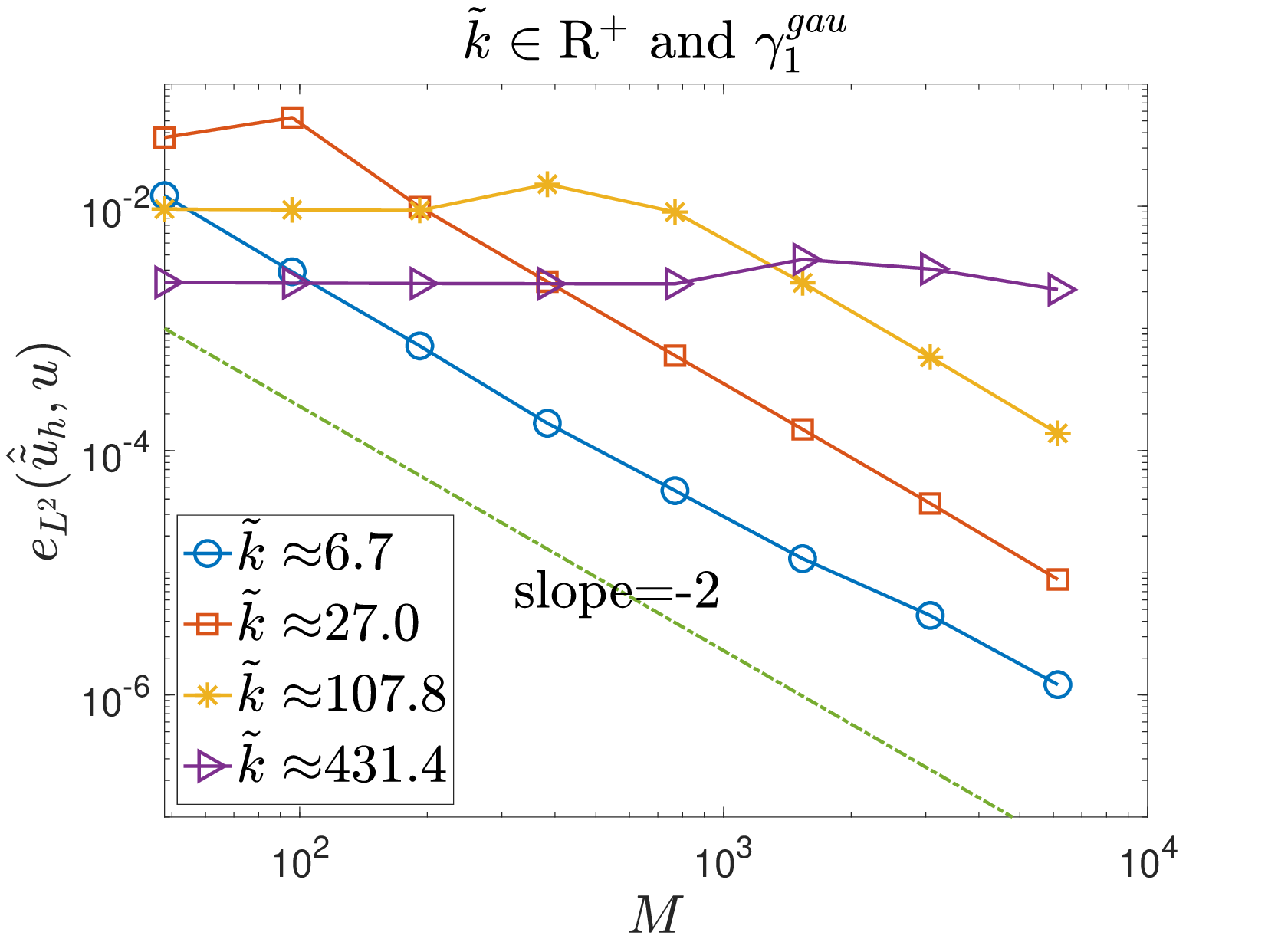}\label{subfig:ex1_fixed_deltaXk_c}}
		\subfigure[]{\includegraphics[width=0.33\textwidth]{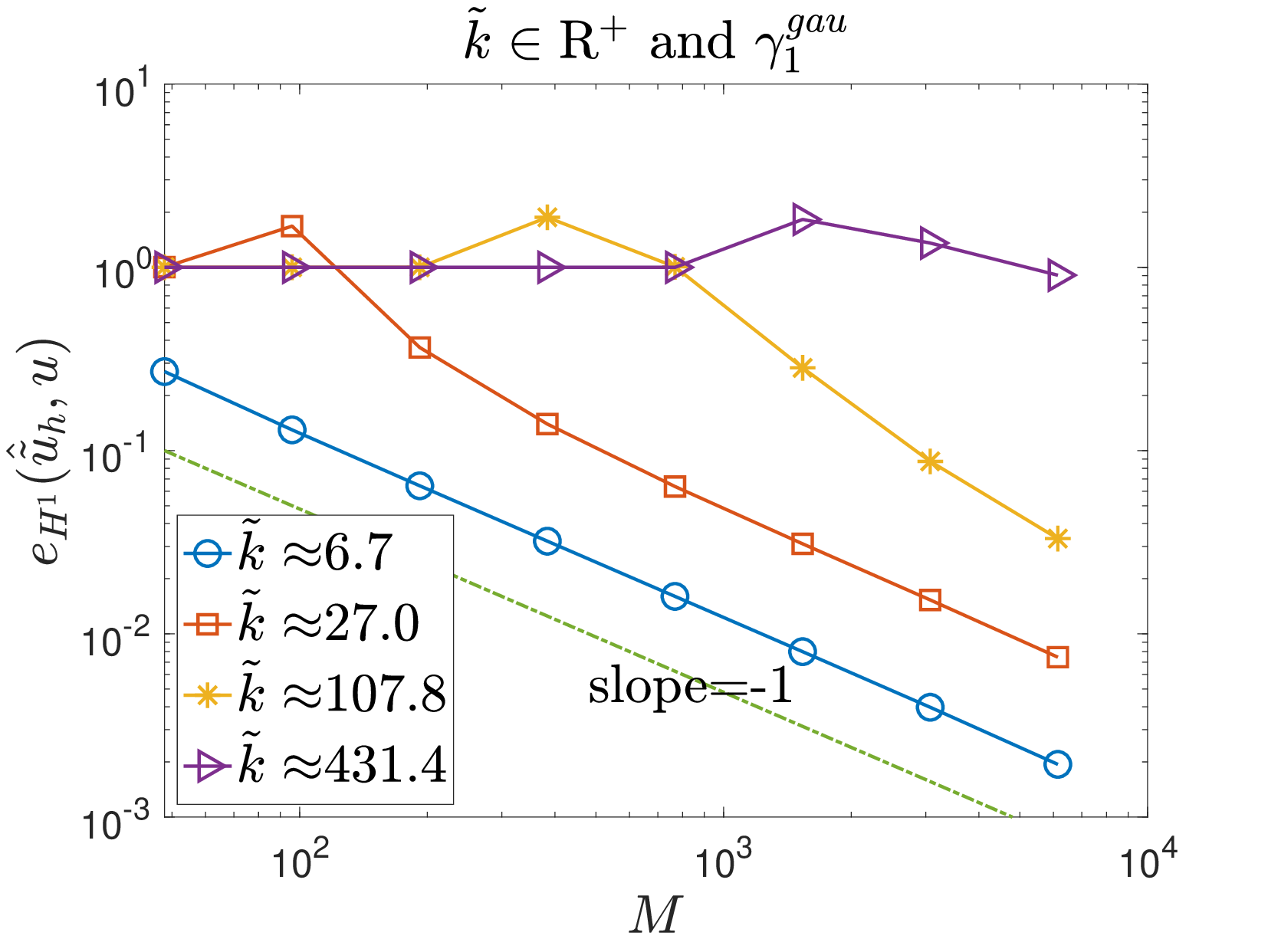}\label{subfig:ex1_fixed_deltaXk_d}}
	\caption{(Example 3) The relative errors $e_{L^2}(\hat{\tilde u}_{h},u)$ and $e_{H^1}(\hat{\tilde u}_{h},u)$ for different $\tilde k\in\R^+$: panels (a)--(b) show errors with the kernel $\gamma_1^{exp}$ and panels (c)--(d) show the errors with the kernel $\gamma_1^{gau}$.} \label{fig:ex1_fixed_deltaXk}
\end{figure}

We first consider the scheme~\eqref{eq:DAScase1eq1}--\eqref{eq:DAScase1eq2} solving nonlocal PML problems with $\tilde k\in \R^+$. Figure~\ref{fig:ex1_fixed_deltaXk} plots the relative errors $e_{L^2}(\hat{\tilde u}_{h},u)$ and $e_{H^1}(\hat{\tilde u}_{h},u)$ for kernels $\gamma_1^{exp}$ and $\gamma_1^{gau}$, respectively.

In subfigures~\ref{subfig:ex1_fixed_deltaXk_a}--\ref{subfig:ex1_fixed_deltaXk_b} the kernel is $\gamma_1^{exp}$ and $(k,\delta)=(\frac{\pi}{5},\frac{1}{2\pi})$, $(\frac{4\pi}{5},\frac{1}{8\pi})$, $(\frac{16\pi}{5},\frac{1}{32\pi})$, $(\frac{64\pi}{5},\frac{1}{128\pi})$, which implies $\tilde k\approx6.3$, $25.3$, $101.0$ and $404.1$ by Eq.~\eqref{eq:ktildekdelta}.

In subfigures~\ref{subfig:ex1_fixed_deltaXk_c}--\ref{subfig:ex1_fixed_deltaXk_d} the kernel is $\gamma_1^{gau}$ and $(k,\delta)=(\frac{\pi}{5},\frac{5}{\pi})$, $(\frac{4\pi}{5},\frac{5}{4\pi})$, $(\frac{16\pi}{5},\frac{5}{16\pi})$, $(\frac{64\pi}{5},\frac{5}{64\pi})$, which implies $\tilde k\approx6.7$, $27.0$, $107.8$ and $431.4$ by Eq.~\eqref{eq:gausskerneltildek}.

\begin{figure} 
	\centering
		\subfigure[]{\includegraphics[width=0.33\textwidth]{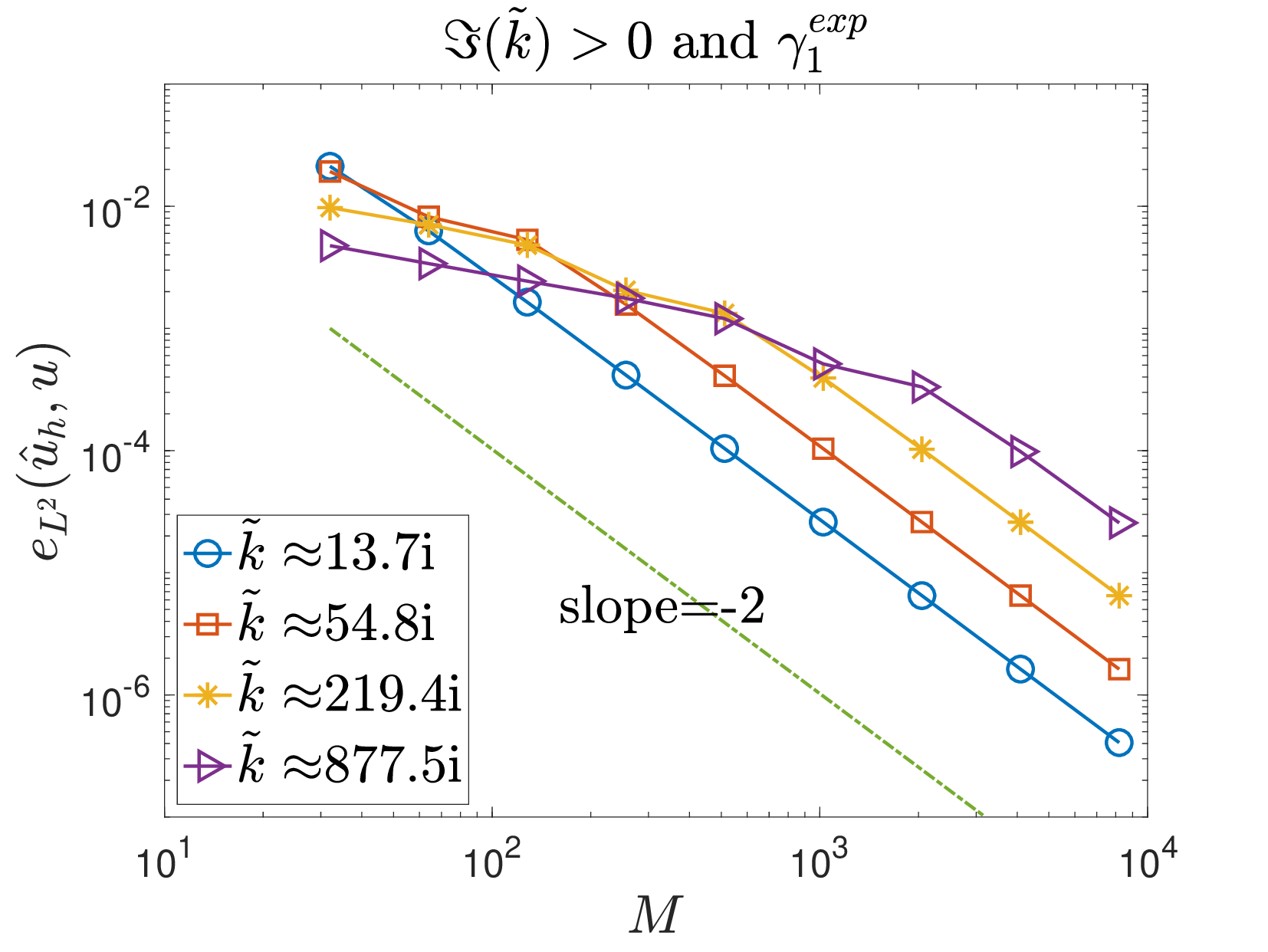}\label{subfig:ex1_fixed_deltaXk_gt1_a}}
		\subfigure[]{\includegraphics[width=0.33\textwidth]{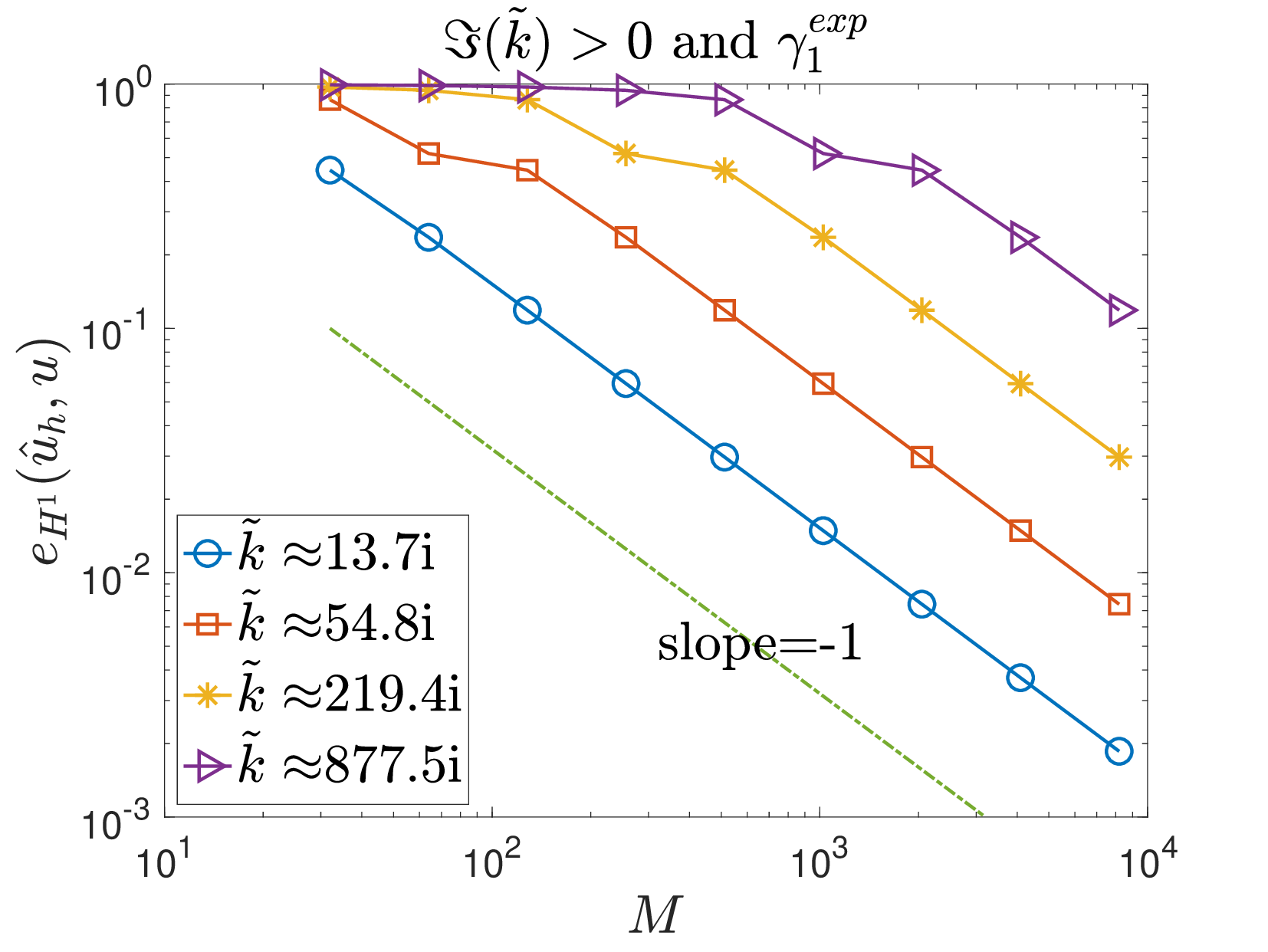}\label{subfig:ex1_fixed_deltaXk_gt1_b}} \\
		\subfigure[]{\includegraphics[width=0.33\textwidth]{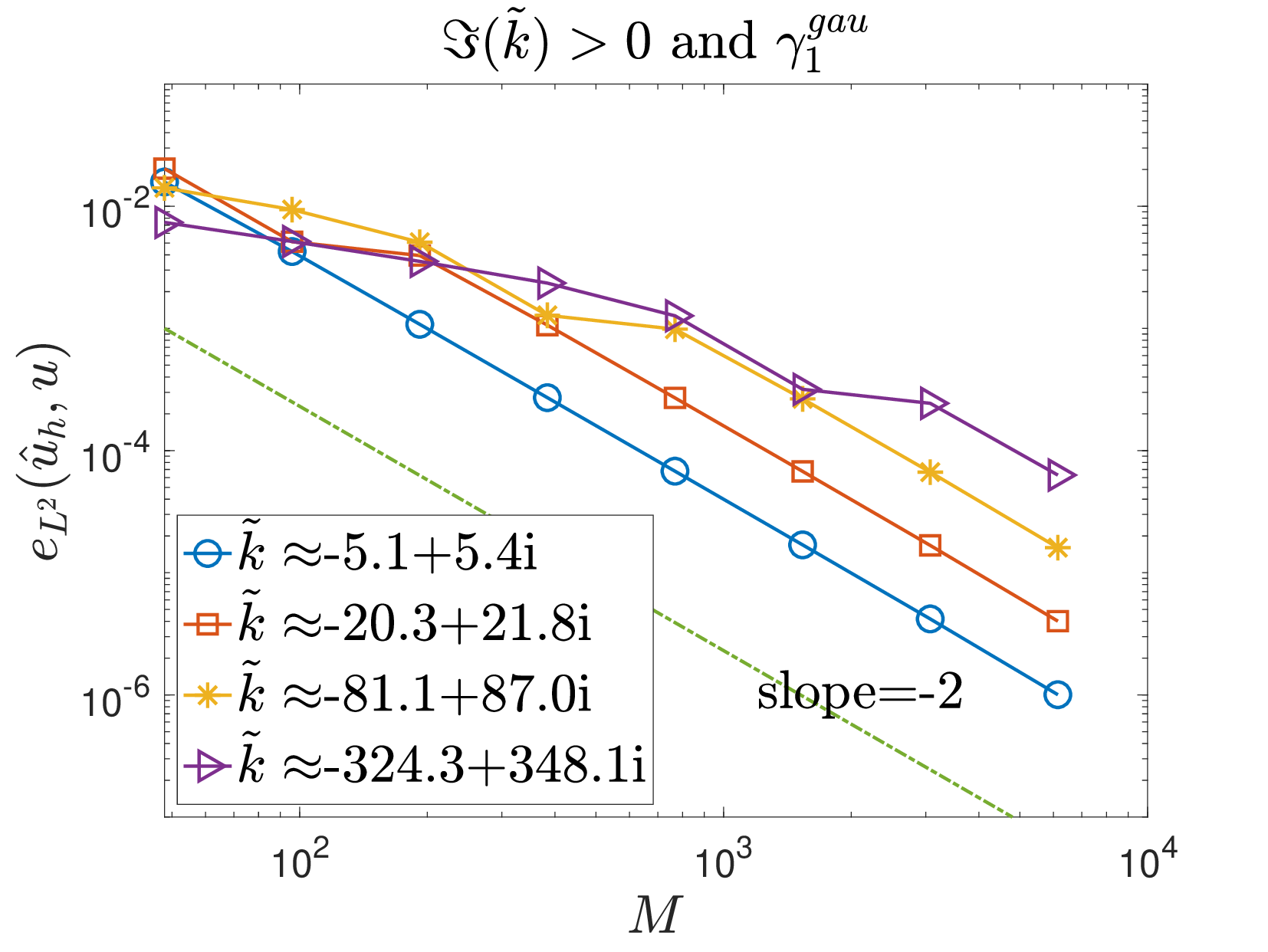}\label{subfig:ex1_fixed_deltaXk_gt1_c}}
		\subfigure[]{\includegraphics[width=0.33\textwidth]{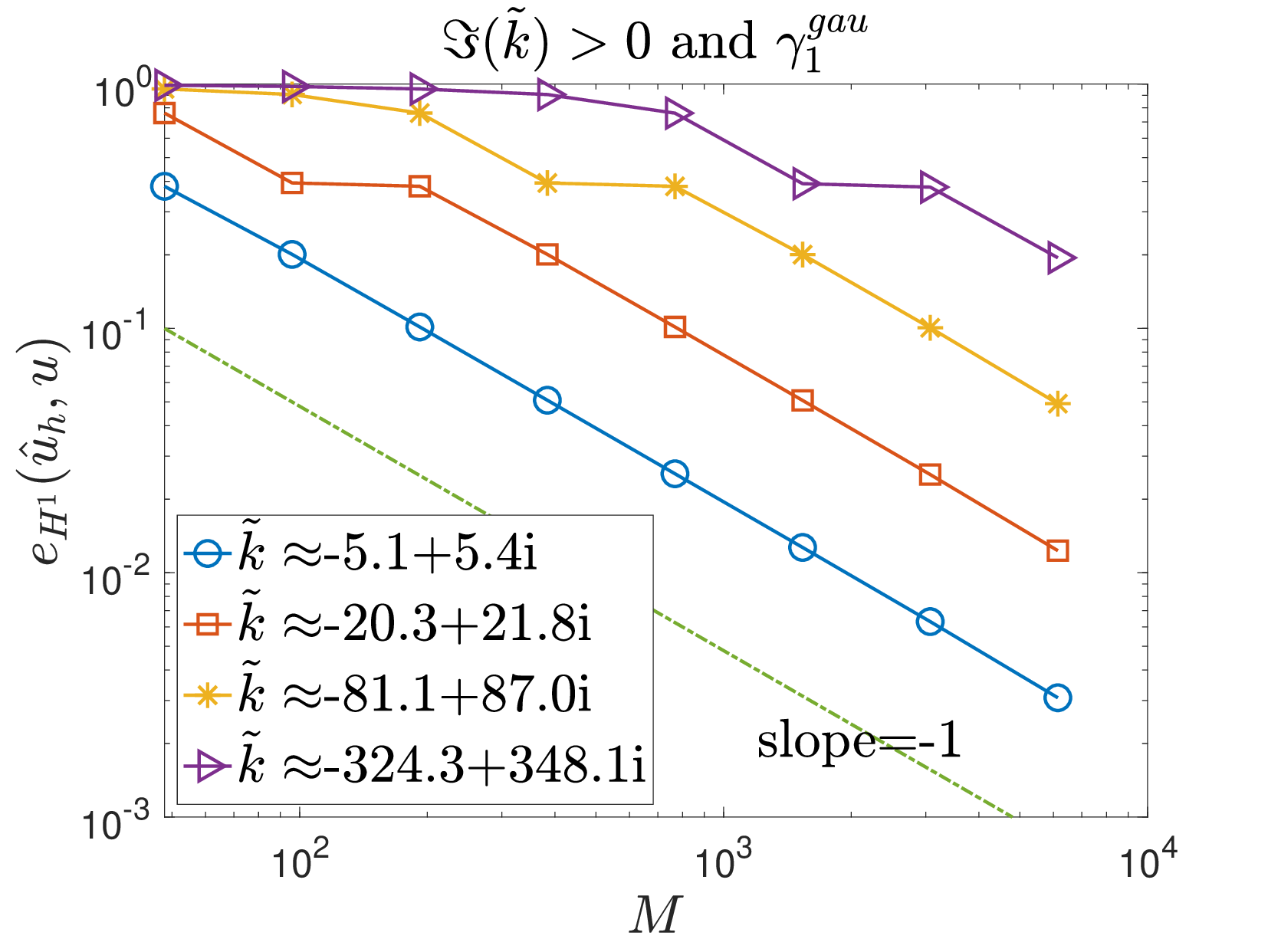}\label{subfig:ex1_fixed_deltaXk_gt1_d}}
	\caption{(Example 3) The relative errors $e_{L^2}(\hat{u}_{h},u)$ and $e_{H^1}(\hat u_{h},u)$ for different $\tilde k$ satisfying $\Im(\tilde k)>0$: panels (a)--(b) show errors with kernel $\gamma_1^{exp}$ and panels (c)--(d) show errors with kernel $\gamma_1^{gau}$.} \label{fig:ex1_fixed_deltaXk_gt1}
\end{figure}

We then consider the scheme~\eqref{eq:DAScase2eq1}--\eqref{eq:DAScase2eq2} solving nonlocal Helmholtz problems with $\Im(\tilde k)>0$. Figure~\ref{fig:ex1_fixed_deltaXk_gt1} plots the relative errors $e_{L^2}(\hat{u}_{h},u)$ and $e_{H^1}(\hat{u}_{h},u)$ for kernels $\gamma_1^{exp}$ and $\gamma_1^{gau}$, respectively.

In subfigures~\ref{subfig:ex1_fixed_deltaXk_gt1_a}--\ref{subfig:ex1_fixed_deltaXk_gt1_b} the kernel is $\gamma_1^{exp}$ and $(k,\delta)=(\frac{\pi}{5},\frac{11}{2\pi})$, $(\frac{4\pi}{5},\frac{11}{8\pi})$, $(\frac{16\pi}{5},\frac{11}{32\pi})$, $(\frac{64\pi}{5},\frac{11}{128\pi})$, which implies $\tilde k\approx13.7\i$, $54.8\i$, $219.4\i$ and $877.5\i$ by Eq.~\eqref{eq:ktildekdelta}.

In subfigures~\ref{subfig:ex1_fixed_deltaXk_gt1_c}--\ref{subfig:ex1_fixed_deltaXk_gt1_d} the kernel is $\gamma_1^{gau}$ and $(k,\delta)=(\frac{\pi}{5},\frac{15}{\pi})$, $(\frac{4\pi}{5},\frac{15}{4\pi})$, $(\frac{16\pi}{5},\frac{15}{16\pi})$, $(\frac{64\pi}{5},\frac{15}{64\pi})$, which implies $\tilde k\approx-5.1+5.4\i$, $-20.3+21.8\i$, $-81.1+87.0\i$ and $-324.3+348.1\i$ by Eq.~\eqref{eq:gausskerneltildek}.

Figures~\ref{fig:ex1_fixed_deltaXk} and \ref{fig:ex1_fixed_deltaXk_gt1} show that for small $k$, the errors in $L^2$-norms and $H^1$-seminorms decay at rates $2$ and $1$ in the log-log scale, respectively; and for larger $k$, these errors first stay around some constants, and converge at the rate $2$ (or $1$) for sufficiently small $h$.

\subsection{Asymptotic compatibility of the discrete scheme}\label{subsec:AC}
\begin{figure} 
	\centering
	\includegraphics[width=0.33\textwidth]{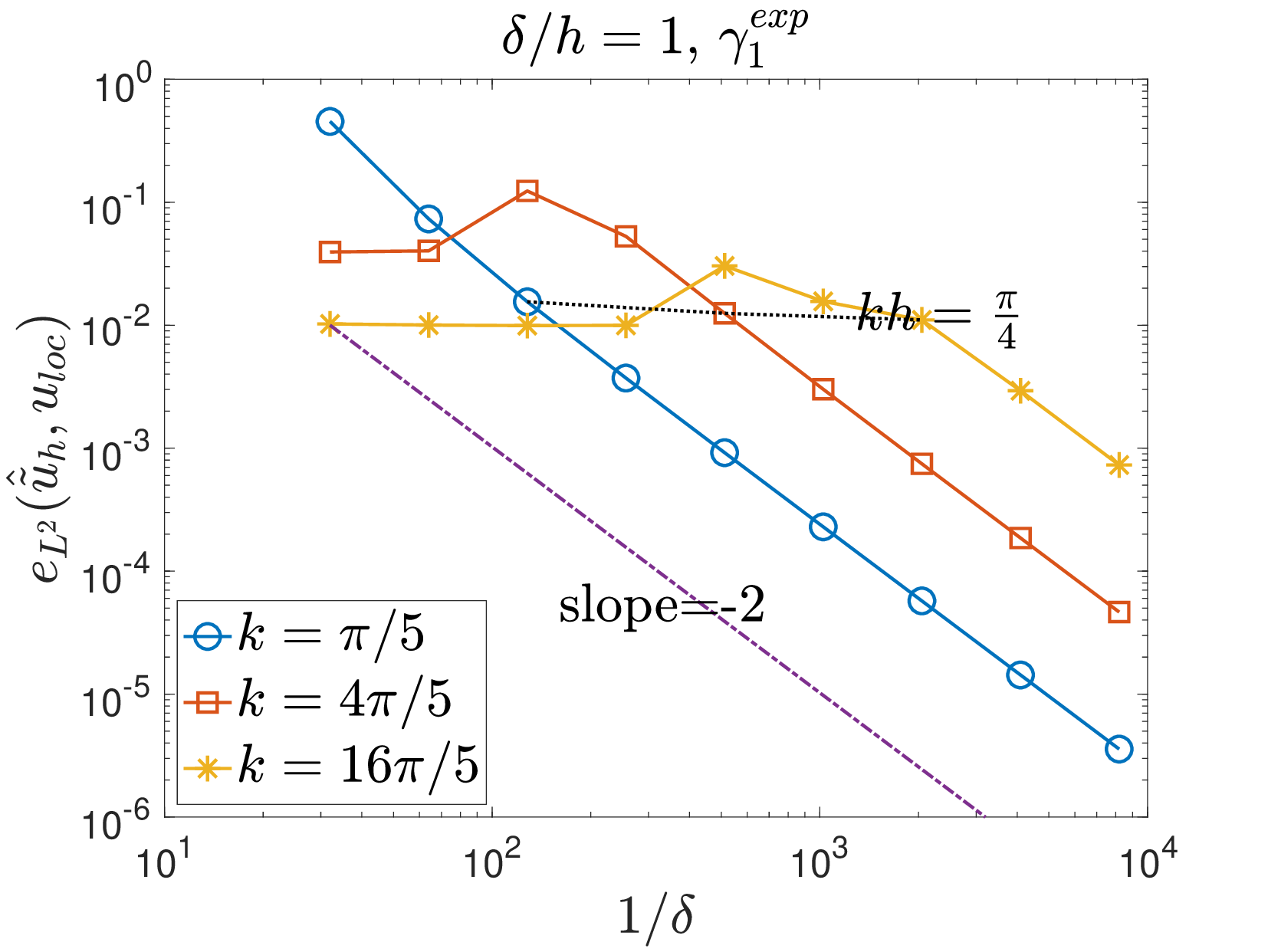}
	\includegraphics[width=0.33\textwidth]{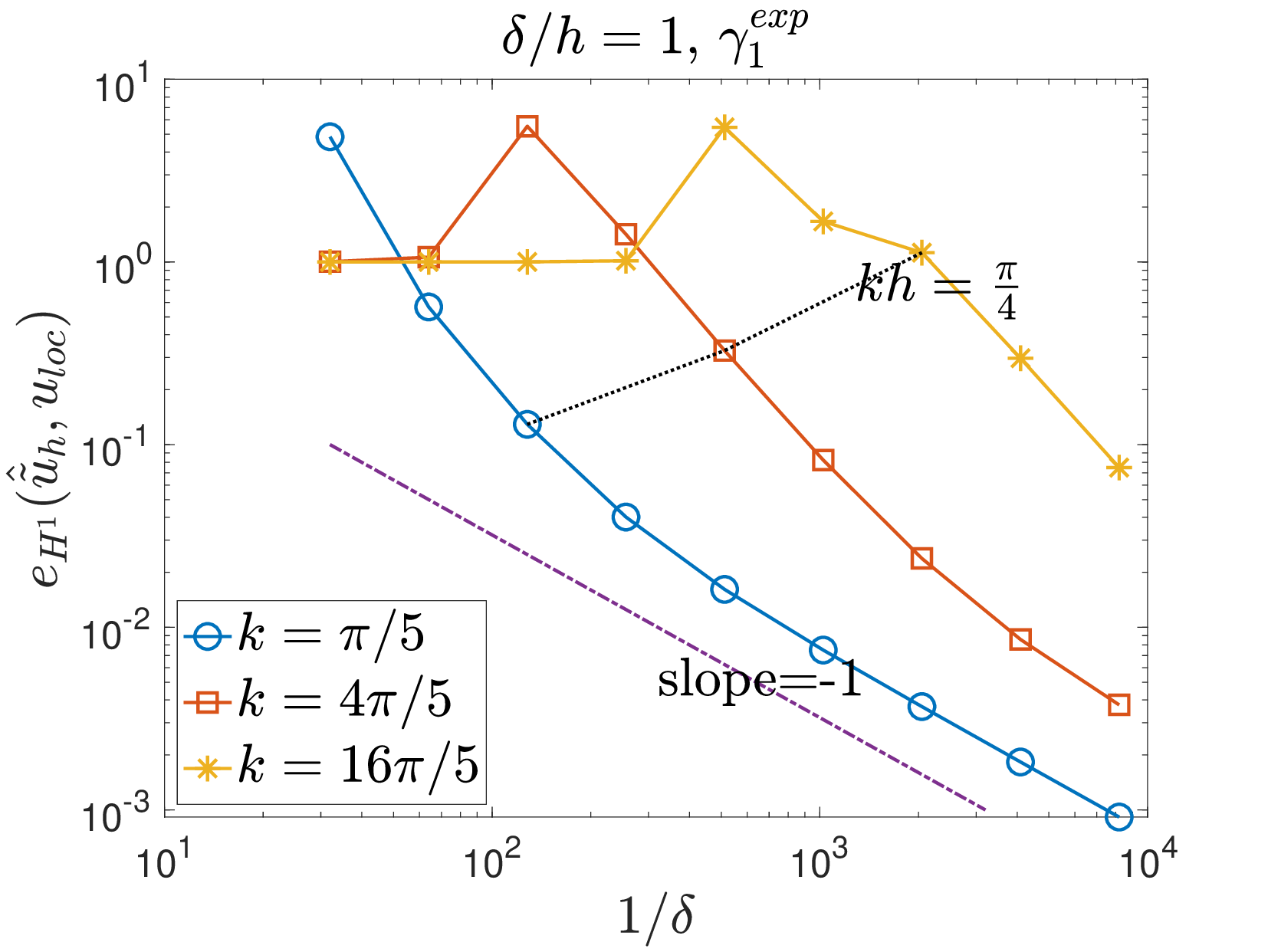}\\
	\includegraphics[width=0.33\textwidth]{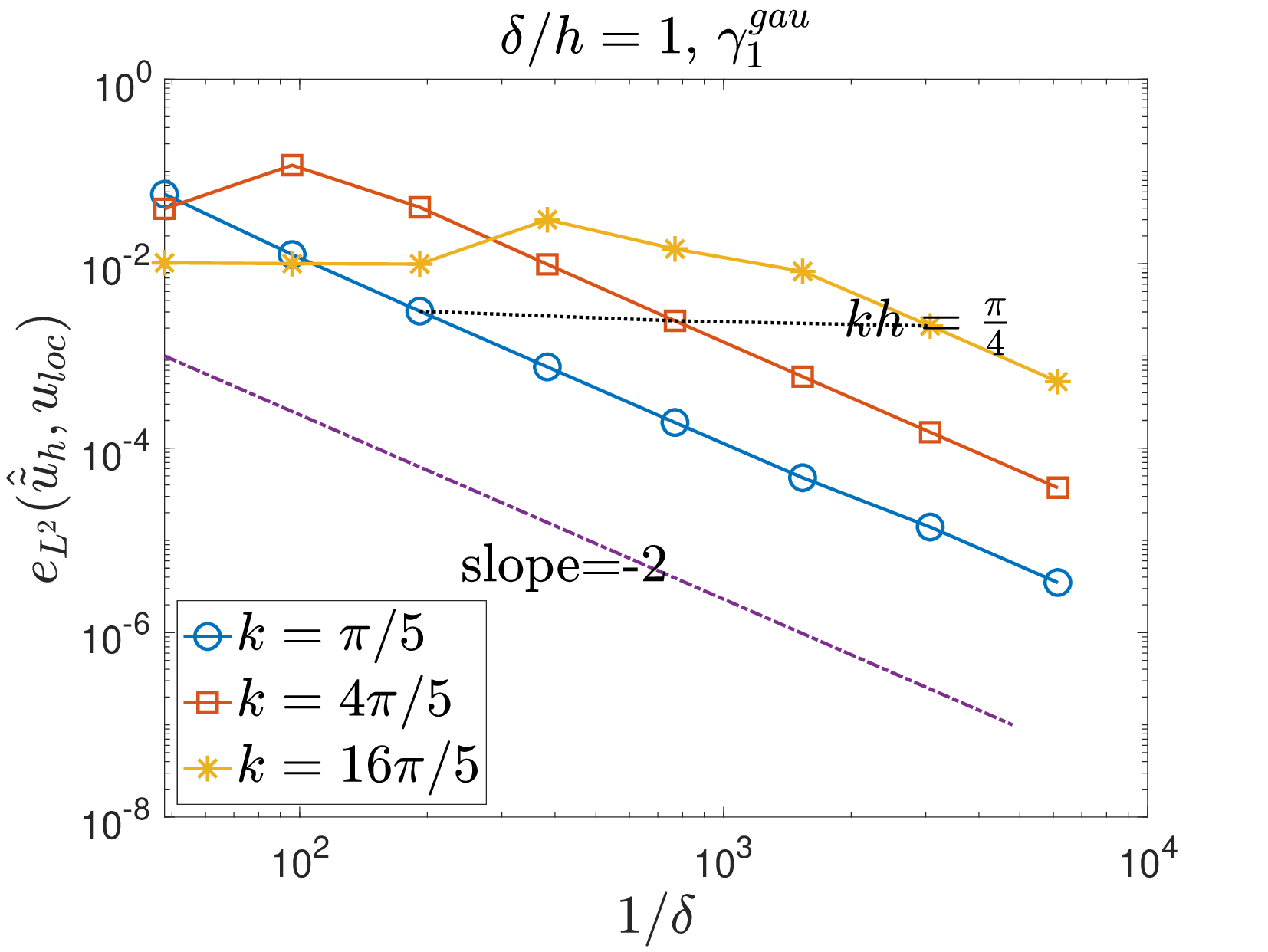}
	\includegraphics[width=0.33\textwidth]{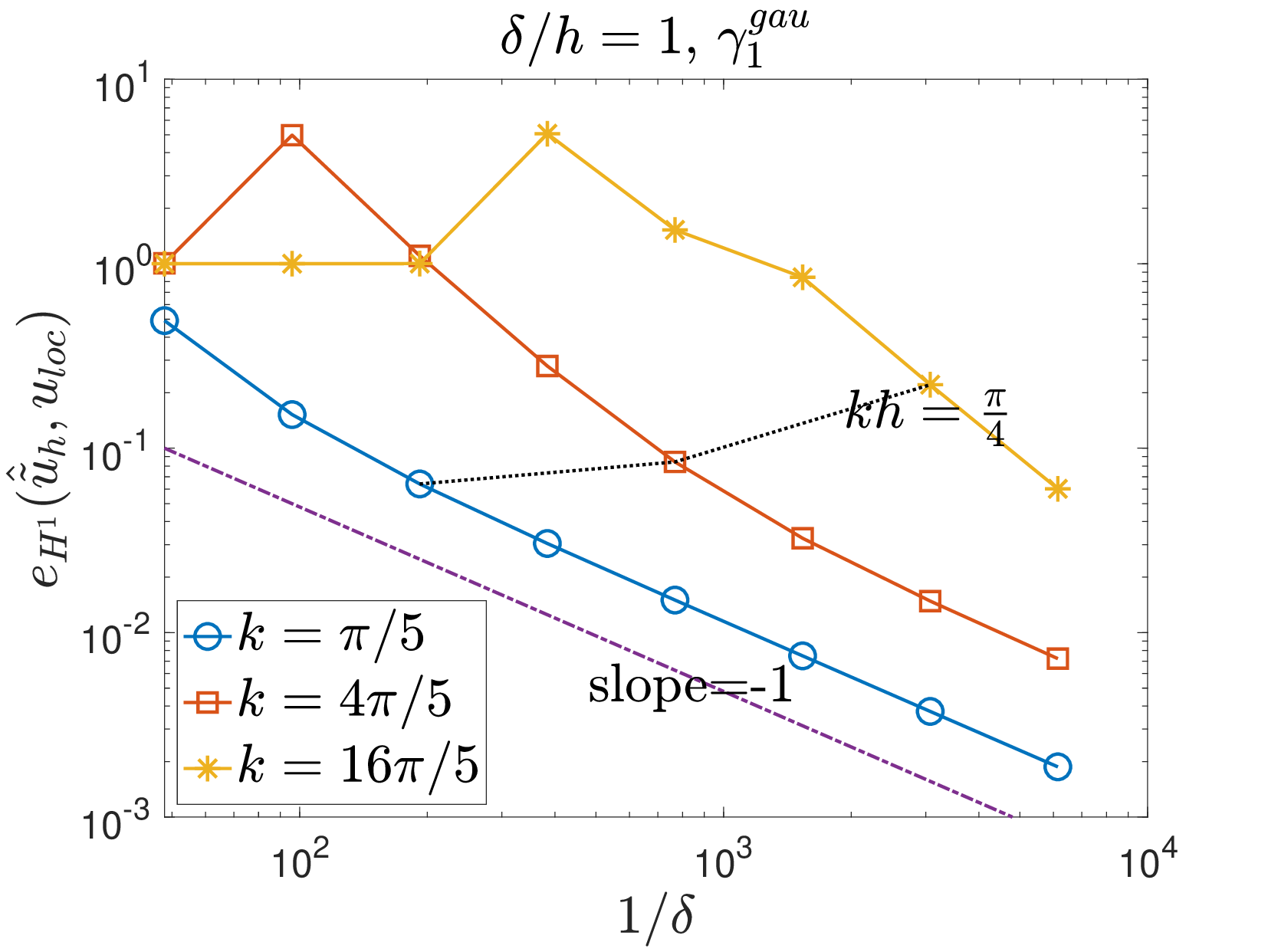}
	\caption{(Example 4) Asymptotically compatibility: the errors $e_{L^2}(\hat{\tilde u}_{h},u_{loc})$ and $e_{H^1}(\hat{\tilde u}_{h},u_{loc})$ between nonlocal numerical solution $\hat{\tilde u}_{h}$ and local solution $u_{loc}$ by fixed $\delta =h$ and taking $\delta \to 0$. 
	} \label{fig:ex1_fixedDelta_heq1}
\end{figure}

{\bf Example 4.} As shown in \eqref{eq:NTL}, the nonlocal operator with PML modifications converges to the corresponding local operator while $\delta \rightarrow 0$. It is interesting to demonstrate whether the numerical solutions $\hat{\tilde u}_{h}$ converge to the corresponding local PML solution $\tilde u_{loc}$ of \eqref{eq:localPML} as both $\delta$ and $h\to0$. To illustrate the quantitative pictures of the consistency of numerical solutions in the local limit, we consider the so-called ``$\delta$-convergence'' by fixing $\delta=h$. In this situation, one has $\delta \rightarrow 0$ while $h\to 0$. Figure~\ref{fig:ex1_fixedDelta_heq1} plots the relative errors of both the $L^2$- and $H^1$-seminorm for different $k=\frac{\pi}{5}$, $\frac{4\pi}{5}$ and $\frac{16\pi}{5}$. One can that the convergence rates are the second order for $L^2$-errors and the first order for $H^1$-errors with respect to $\delta$, which is consistent to the convergence analysis in \cite{TianDu}. 
 

\section{Conclusion and discussion}
In this paper the PML equation is derived from the weak form of the nonlocal Helmholtz equation. In terms of the weighed average values, we theoretically prove the decay properties of the nonlocal Helmholtz solution and the nonlocal PML solution with a general kernel.  For a typical kernel $\gamma_1(s)=\frac12 e^{-| s|}$, we may formulate the Green's function for the nonlocal Helmholtz equation, and give a refined stability estimate on the decay properties of solutions. The exact formula of Green's function not only plays an important role on the stability estimates of the nonlocal solutions, but also shows an intimate connection between local and nonlocal solutions, i.e., they can be expressed by each other. Numerical examples are given to demonstrate the effectiveness of our PLM strategy and to verify the theoretical findings. As far as we know, this is a pioneering work to provide a general methodology for systematically constructing the PML for the nonlocal Helmholtz and analyzing the corresponding solution properties.

Although our PML technique is derived and the behaviors of solutions are studied in homogeneous medium, it also can be applied to the nonlocal Helmholtz equation in multi-scale mediums. In the future, we will extend to our PML strategy to 2D and 3D cases and other nonlocal wave propagations such as nonlocal wave equations.

\section*{Appendix: The proof  of Proposition \ref{Le1}.}

As the well-documented dispersion analysis and numerical tests \cite{Du2016handbook,du2018nonlocal}, the nonlocal wave propagation is different from the local wave propagation. We assume that the nonlocal wave is $e^{\i\tilde k x}$ for some $\tilde k\in\mathbb{C}$ and large $x\in\R^+$ when $f\equiv0$, which is different from the local wave $e^{\i kx}$. $\tilde k$ should make $e^{\i\tilde k x}$ satisfy the nonlocal Helmholtz equation, that is
\begin{align}
\mathcal{L}_\delta e^{\i\tilde kx} - k^2 e^{\i \tilde kx} =0,
\end{align}
which implies that $\tilde k$ is the solution to the equation $\mu(\tilde k)=k^2$ (see \eqref{eq:relation}). We choose $\tilde k\in\R^+$ or $\Im(\tilde k)\in\R^+$ in the set of solutions to \eqref{eq:relation}. The reason is that when $\tilde k\in\R^+$, $e^{\i\tilde k x}$ is an outgoing wave and when $\Im(\tilde k)\neq0$, $\Im(\tilde k)\in\R^+$ makes $e^{\i\tilde k x}$ be a damped wave. These phenomena have been studied and can be observed in \cite{Du2016handbook,du2018nonlocal}. In the literature, the Gaussian kernel was studied and it's shown that for relatively small $k$ corresponding to $\tilde k\in\R^+$ here, the group velocity of the nonlocal wave is almost equal to the phase velocity, and thus the nonlocal wave travels like the sine wave. However, for relatively large $k$ corresponding to $\Im(\tilde k)\in\R^+$, the group velocity is close to zero, and as a result, the amplitude of the nonlocal wave approaches to zero as $x$ goes to infinity.

By using $\tilde k$ introduced above, we define an auxiliary function $G_{x_0}(x)$ (see \eqref{eq:green}). Then we multiply Eq.~\eqref{eq:nonlocalHelmholtz} by $G_{x_0}(x)$ with a direct calculation, and obtain 
\begin{align*}
	 & \int_{\R} G_{x_0}(x)f(x)\dx = \int_{\R} G_{x_0}(x) \big(\mathcal{L}_\delta u(x)-k^2u(x)\big) \dx                                \\
= & \int_{\R} G_{x_0}(x)u(x) \left[ \frac{1}{\delta^2}\int_{\R} \bigg(1-e^{\i\tilde k\delta s}\bigg) \gamma_1(s) \ds -k^2 \right] \dx                          \\
	 & - \frac{1}{\delta^2} \int_{-\infty}^{x_0} u(x) \left[ \int_{-\infty}^{\frac{x-x_0}{\delta}} \left( C_2 e^{\i\tilde{k}(x-\delta s)} - C_1 e^{-\i\tilde{k}(x-\delta s)} \right) \gamma_1(s) \ds \right]\dx         \\
	 & - \frac{1}{\delta^2} \int_{x_0}^{+\infty} u(x) \left[ \int_{\frac{x-x_0}{\delta}}^{+\infty} \left( C_1 e^{-\i\tilde{k}(x-\delta s)} - C_2 e^{\i\tilde{k}(x-\delta s)} \right) \gamma_1(s) \ds \right]\dx. 
\end{align*}
Using the relation \eqref{eq:relation} and taking $t= x-x_0$, we further have 
\begin{align*}
	 & \int_{\R} G_{x_0}(x)f(x)\dx \\
	= & \frac{1}{\delta^2} \int_{-\infty}^{0} u(t+x_0) \left[ \int_{-\infty}^{\frac{t}{\delta}} \frac{1}{2\i\tilde k} \left(e^{\i\tilde k(t - \delta s)} - e^{-\i\tilde k(t-\delta s)} \right) \gamma_1(s) \ds \right]\mathrm{d}t \\
	 &+ \frac{1}{\delta^2} \int_{0}^{+\infty} u(t+x_0) \left[\int_{\frac{t}{\delta}}^{+\infty} \frac{1}{2\i\tilde k} \left( e^{-\i\tilde k (t-\delta s)} - e^{\i\tilde k(t-\delta s)} \right) \gamma_1(s) \ds \right]\mathrm{d}t \\
	=: & \int_{-\infty}^{0} u(t+x_0) w_1(t) \mathrm{d}t + \int_{0}^{+\infty} u(t+x_0) w_2(t) \mathrm{d}t ,
\end{align*}
where where $w_1(t)$ and $w_2(t)$ are given in \eqref{w1} and \eqref{w2}. It is direct to verify $w_1(-t)=w_2(t)$ for $t>0$. Set $\kappa(t) = w_1(t)$ for $t<0$ and $\kappa(t) = w_2(t)$ for $t>0$. Thus, we have $\kappa(t) = \kappa(-t)$ and 
\begin{align} \label{eq:barqeq} 
	 & \int_{\R} u(t+x_0) \kappa(t) \mathrm{d}t = \int_{\R} G_{x_0}(x)f(x)\dx.
\end{align}

Noting the definition of the weighted average value $u^w(x)$ in \eqref{eq:barqdef}, and using relation \eqref{eq:barqeq} for any point $x_0$, we immediately have the identity \eqref{eq:baruExp}. 

\bibliographystyle{siam}

\end{document}